\documentclass[a4paper]{scrartcl}
       
\usepackage[utf8]{inputenc}
\usepackage[T1]{fontenc}
\usepackage[english]{babel} 
\usepackage{lmodern}
\usepackage[onehalfspacing]{setspace}
\usepackage{amsmath,amssymb,amsthm,amsfonts,amsbsy,latexsym}
\usepackage{graphicx}
\usepackage{natbib}
\usepackage{geometry}
\usepackage{float}
\usepackage[]{algorithm2e}
\usepackage{pgf,tikz}
\usepackage[pdfpagelabels=true]{hyperref}

\newcommand{\E}{\mbox{I\negthinspace E}}

\newcommand{\R}{\mathbb{R}}

\DeclareMathOperator*{\argmin}{argmin}

\DeclareMathOperator*{\sign}{sign}

\DeclareMathOperator*{\rk}{rk}

\providecommand{\keywords}[1]
{
  \small	
  \textbf{\textit{Keywords---}} #1
}

\numberwithin{equation}{section}

\newtheorem{thm}{Theorem}[section]
\newtheorem{Def}{Definition}[section]
\newtheorem{lem}{Lemma}[section]
\newtheorem{rem}{Remark}[section]
\newtheorem{cor}{Corollary}[section]
\newtheorem{ex}{Example}[section]
\newtheorem{ass}{Assumption}[section]
\newtheorem*{bew}{Proof}

\begin{document}

\author{Tino Werner}
\date{\today}
\title{Quantitative robustness of instance ranking problems} 
\maketitle

\begin{abstract} Instance ranking problems intend to recover the true ordering of the instances in a data set with a variety of applications in for example scientific, social and financial contexts. Robust statistics studies the behaviour of estimators in the presence of perturbations of the data resp. the underlying distribution and provides different concepts to characterize local and global robustness. In this work, we concentrate on the global robustness of parametric ranking problems in terms of the breakdown point which measures the fraction of samples that need to be perturbed in order to let the estimator take unreasonable values. However, existing breakdown point notions do not cover ranking problems so far. We propose to define a breakdown of the estimator as a sign-reversal of all components which causes the predicted ranking to be potentially completely inverted, therefore we call our concept the order-inversal breakdown point (OIBDP). We will study the OIBDP, based on a linear model, for several different ranking problems that we carefully distinguish and provide least favorable outlier configurations, characterizations of the order-inversal breakdown point as well as sharp asymptotic upper bounds. We also outline the case of SVM-type ranking estimators. \end{abstract}

\keywords{Breakdown point; Quantitative robustness; Instance ranking problems; Sparsity} 

\section{Introduction} 

A well-known issue when analyzing data is that the data usually is not clean but consists of perturbations that can severely distort an estimator. Instances that are distant from the majority of the data are often termed as ''outliers''. For a well-founded analysis, it is required to define an underlying ideal model so that the data points are interpreted as independent realizations from this model. However, wrong model assumptions let the real data appear as contaminated data. In these cases, methods from robust statistics (\cite{huber}, \cite{hampel}, \cite{rieder}, \cite{maronna}) to handle these phenomena are necessary to incorporate even contaminated data points appropriately since just removing possible outliers is the wrong way as for example discussed in \cite{hampel}. Even worse, in contrast to the classical convex contamination model (see e.g. \cite[Sec. 4.2]{rieder}) where an instance either stems from a contaminated distribution or from the ideal distribution, the cell-wise outlier model from \cite{alqallaf} allows for contaminating the single predictor components for each instance independently which causes the probabiliity to have even one clean instance in the data to tend to zero which again is a manifestation of the curse of dimensionality. 

Robust statistics provides two concepts to measure the quantitative robustness of an estimator. Since robust statistics identifies estimators as statistical functionals (\cite{huber}, \cite{hampel}, \cite{rieder}, \cite{maronna}), functional derivatives (e.g., \cite{averbukh}) can be applied in order to linearize this functional in a first-order expansion which goes back to \cite{mises}. The functional derivative, usually a G\^{a}teaux derivative, has been identified in \cite{hampel74} with the influence curve which is an important diagnostic tool which measures the infinitesimal impact of one data point on the estimator. %and which also has been considered for example in \cite{reeds}, \cite{clarke} and \cite{rieder}. The concept of influence curves is well-known, but its success is based on the much deeper theoretical fact that standard M-estimators as well as for example MD- and R-estimators allow for an asymptotic linear expansion (see e.g. \cite{fernholz}). 

In contrast to the influence curve which quantifies the local robustness of an estimator, i.e., only allowing for an infinitesimal fraction of the data being contaminated, the breakdown point (BDP), introduced in \cite[Sec. 6]{hampel71} in a functional version and in \cite{huber83} in a finite-sample version, studies the global robustness of an estimator. The finite-sample BDP from \cite{huber83} quantifies the minimum fraction of instances in a data set so that contaminating any such fraction of data points arbitrarily can let the estimator ''break down'' while the functional BDP essentially quantifies the allowed maximum Prokhorov distance between the ideal and the contaminated distribution without the estimator breaking down. There has yet been a lot of work on BDPs, see for example \cite{rous84}, \cite{rous85}, \cite{davies93} and \cite{hubert97}, \cite{genton98}, \cite{becker99}, \cite{gather} or \cite{donoho06} which cover location, scale, regression and spatial estimators and \cite{hubert08} who study the BDP for multivariate estimators. Recently, a BDP for classification (\cite{zhao18}) and for multiclass-classification (\cite{qian}) have been proposed.

While regression and classification aim for an exact fit of the response value, there exist types of problems where one is only interested in an ordering of the instances and not of the particular responses. These problems are ranking problems which are very important in for example in document ranking (\cite{page}, \cite{herb}, \cite{cao}), medicine (\cite{agar09}), credit risk-screening (\cite{clem13}) or biology and chemistry (\cite{agar10}, \cite{kayala}, \cite{kayala11}, \cite{morrison}). Due to the global nature of ranking problems where essentially each instance pair is compared, the existing global robustness measures, i.e., the existing BDP concepts, are not suitable here. 

Consider the problem to order instances in a data set. If responses are available, this can be identified with minimizing a pair-wise loss function, i.e., which operates on pairs of responses and their predictions by checking if their ordering coindices, as shown in the seminal work of \cite{herb}, \cite{herb99}. We then speak of instance ranking problems in the terminology of \cite{furn11}, in contrast to object or item ranking problems where no responses are available (e.g., \cite{cohen99}, \cite{huller15}) and label ranking problems where a multicategorical response is given and the goal is for each instance to recover the ordering of the probabilities that the response belongs to the respective class (\cite{furn08}, \cite{huller10b}). \cite{clem08} proposed the statistical framework for such instance ranking problems which emerge from ordinal regression (\cite{herb}) and proved that the common approach of empirical risk minimization (ERM) is indeed suitable for such ranking problems. There are three ways of casting a ranking problem, i.e., either the ordering of all instances has to be correct (hard ranking), one just wants to identify the top $K$ instances for a given $K$ (weak ranking, \cite{clem08b}) or the best $K$ instances have to be identified and the ordering of these instances has to be correct (localized ranking, \cite{clem08b}). Furthermore, one distinguishes between binary responses which lead to binary or bipartite ranking problems which have been considered in many works (e.g., \cite{joachims02}, \cite{freund}, \cite{clem08d}), $d-$partite ranking problems for categorical responses with $d$ categories (e.g., \cite{clem12}, \cite{furn09}) and continuous ranking problems (\cite{sculley}, \cite{clem18}).  

Instance ranking problems are usually solved by learning a real-valued, here parametric, scoring function which assigns a score to each instance with the goal to minimize some ranking error between the predicted ordering of the instances according to the scores and the true ordering. The peculiarity of ranking problems is that they have an inherent equivariance nature, i.e., multiplying each response with the same positive factor or adding the same fixed value to it does not alter the ordering. The regression BDP is essentially understood in the sense that at least such a high fraction of outliers has full control over  the parameter, more precisely, any finite bound on the norm can be exceeded, though achieving an infinite norm is in general not possible. It is out of question that for a ranking prediction, it would be even worse to predict an inverted ordering than to perform random guessing (to which infinite coefficients essentially correspond) which exactly motivates our so-called OIBDP which is the minimum fraction of perturbed data points so that the non-zero coefficient components can be inverted. At the first glance, we get unreasonable BDPs in high-dimensional settings, i.e., if the predictor dimension is no longer smaller than the number of observations, but this can be remedied by assuming sparse underlying models resp. sparse model selection which is natural in such settings (e.g., \cite{bu}). 

Our contribution is threefold: \textbf{i)} We propose the definition of the order-inversal BDP for ranking problems which embeds the BDP concept of robust statistics into that area of machine learning; \textbf{ii)} we provide explicit worst-case outlier configurations and \textbf{iii)} we compute upper bounds for the corresponding OIBDPs for different ranking problems. 

The rest of this work is organized as follows. Sec. \ref{prelim} compiles necessary preliminaries in terms of a more concise definition of the loss functions corresponding to the different ranking problems as well as the BDP concept. In Sec. \ref{motivsec}, we show why neither the classical BDP for regression nor the angular BDP for classification are suitable for ranking problems and propose the OIBDP for ranking. In Sections \ref{hardcontsec}, \ref{hardbinsec} and \ref{locsec}, we propose outlier schemes and prove asymptotic bounds for the OIBDP for hard continuous resp. hard binary and hard $d-$partite resp. localized continuous ranking problems. In Sec. \ref{othersec}, we discuss the applicability of BDP concepts to the remaining instance ranking problems. In Sec. \ref{discusssec}, we relate the computed BDPs to sparse underlying models and outline how robust ranking can be achieved in practical applications. Sec. \ref{furthersec} is an outlook devoted to SVM- resp. SVR-type approaches. Further results and selected proofs can be found in the Appendix.

\section{Preliminaries} \label{prelim} 

We start by revisiting suitable loss functions for different types of instance ranking problems and the breakdown point concept. 

\subsection{Ranking problems} \label{rankingsec} 

Let $D=(X,Y)$ be a data set with regressors $X_i \in \mathcal{X} \subset \R^p$ and responses $Y_i \in \mathcal{Y} \subset \R$. Instance ranking problems, in the following just referred to as ranking problems where no confusion is possible, generally aim at predicting the correct ordering of the instances by predicting the correct ordering of the responses, i.e., $X_i$ will be ranked higher than $X_j$ if $\hat Y_i>\hat Y_j$. Let a scoring function $s_{\beta}: \mathcal{X} \rightarrow \R$ for some parameter $\beta \in \Theta \subset \R^p$ be given. $X_{i,j}$ refers to the $j-$th entry of row $X_i$ and $X_{\cdot,j}$ to the $j-$th column of $X$. 

In the case of hard ranking problems the resulting parametric optimization problem is given by \begin{equation} \label{rankparopt} \min_{\theta \in \Theta}\left(L_n^{hard}(\beta)=\frac{1}{n(n-1)}\mathop{\sum \sum}_{i \ne j} L(X_i,X_j,Y_i,Y_j,\beta) \right) \end{equation} where $L: \mathcal{X} \times \mathcal{X} \times \mathcal{Y} \times \mathcal{Y} \times \Theta \rightarrow [0,\infty]$ is some loss function that compares instance pairs. In \cite{herb99} or \cite{clem08}, $L$ is the indicator function \begin{center} $ \displaystyle L(X_i,X_j,Y_i,Y_j,\beta)=I((Y_i-Y_j)(s_{\beta}(X_i)-s_{\beta}(X_j))<0) $ \end{center} which just checks whether a misranking occurred, i.e., if the true resp. the predicted pair-wise orderings did not coincide, but the actual magnitude of the product is not taken into account. Since this loss function is not even continuous, one often considers surrogate losses (see \cite{TW19b} for an overview). In the following BDP computations, we will always consider general loss functions that can be rewritten as $L((Y_i-Y_j)(s_{\beta}(X_i)-s_{\beta}(X_j))$ in the same manner as classification loss functions are frequently rewritten as $L(ys_{\beta}(x))$.  

For weak ranking problems (\cite{clem08b}), the empirical counterpart of the misclassification risk can be expressed by  \begin{center}  $ \displaystyle L_n^{weak, K}(\beta)=\frac{2}{n}\sum_{i \in Best_K} I(\rk(s_{\beta}(X_i))>K) $  \end{center} with the set $Best_K$ of the true top $K$ indices where the ranks correspond to a descending ordering. Again, the indicator function may be replaced by any classification loss function $L: \mathcal{X} \times \mathcal{Y} \times \beta \rightarrow [0,\infty]$. 
                 
A suitable loss function for localized ranking problems (\cite{clem08b}) is \begin{equation} \label{locrankrisk} L_n^{loc, K}(\beta):=\frac{n-K}{n} L_n^{weak, K}(\beta)+\frac{2}{n(n-1)}\mathop{\sum \sum}_{i<j, i, j \in Best_K} I((s_{\beta}(X_i)-s_{\beta}(X_j))(Y_i-Y_j)<0) .  \end{equation} One can rewrite the second sum by \begin{center} $ \displaystyle \frac{2}{n(n-1)}\mathop{\sum \sum}_{i<j, i, j \in Best_K} I((s_{\beta}(X_i)-s_{\beta}(X_j))(Y_i-Y_j)<0) $ \end{center} or, again, replace the indicator function by surrogates. Note that one may replace the set $\widehat{Best_K}$ in the double sum by $Best_K$. We will discuss both cases in Sec. \ref{locsec}. 

For further discussions of these loss functions and for instance ranking, see \cite{TW19b}.

\subsection{Quantitative robustness} \label{quantrob} 

Robustness of an estimator can be understood in the sense that it allows for perturbations or even large contaminations of the underlying sample without the quality of the estimator being significantly affected. One can distinguish between quantitative and qualitative robustness. The latter goes back to \cite{hampel71} and essentially indicates the continuity of the underlying statistical functionals. \ \\

As for quantitative robustness, one further has to distinguish between global and local robustness. Local robustness is devoted to the effect of small perturbations of the data where the term ''small'' means that, for finite samples, only one observation may be contaminated, so in other words, the influence curve or influence function which is the diagnostic tool for local robustness measures the infinitesimal influence of a single observation on the estimator. In contrast, global robustness allows for large perturbations, i.e., a considerable fraction of the data points being contaminated arbitrarily. The maximum fraction which an estimator can cope with, i.e., without taking unreasonable values, is measured by the breakdown point.

\subsubsection{The breakdown point concept} 

Let $Z_n$ be a sample $(X_1,Y_1),...,(X_n,Y_n)$. Let $\hat \beta(Z_n)$ be the estimated coefficient for the scoring function $s_{\beta}$ based on $Z_n$. The finite-sample BDP of \cite{huber83} is defined as follows. 

\begin{Def} The \textbf{finite-sample breakdown point} of an estimator $\hat \beta$ is defined as \begin{equation} \label{fsbdp} \epsilon^*(\hat \beta,Z_n)=\min\left\{\frac{m}{n} \ \bigg| \ sup_{Z_n^m}(||\hat \beta(Z_n^m)||)=\infty \right\} \end{equation}  where $Z_n^m$ denotes any sample that has exactly $(n-m)$ instances in common with $Z_n$, i.e., $m$ instances can be modified arbitrarily. \end{Def} 

Note that this definition assumes that $\beta \in \R^p$. In cases where $\beta \in \Theta \subset \subset \R^p$, the situation would get more difficult since here a breakdown may be defined in the sense that $\hat \beta$ is located at the boundary of $\Theta$. In this case, one would require some transformation that moves the boundaries of $\Theta$ to infinite values, see for example \cite{he05} who proposed to use a log-transformation for computing the BDP of scale estimators in order to move the value 0 to $-\infty$. %There also exist an asymptotic variant of the finite-sample BDP given in \cite{hampel} which does not consider finite samples but a sequence of estimators based on samples of diverging dimension.
\ \\

A variety of extensions of the BDP concept have been proposed in the literature. \cite{stromberg} and \cite{sakata} proposed BDPs for regression, \cite{sakata98} suggested a BDP definition for location-scale estimators while \cite{genton98} propose the spatial BDP for variogram estimators and \cite{genton03} and \cite{genton03b} introduce a BDP for dependent samples (time series). \cite{donoho06}, \cite{donoho06b} propose a BDP for model selection, \cite{kanamori} study the BDP for SVMs and \cite{hennig08} transferred the BDP concept to the dissolution point concept for clustering. \cite{horbenko12} suggest an expected BDP that respects the ideal distribution of the original data. See \cite{davies05} for a notable discussion paper on BDPs.

\subsubsection{Angular breakdown point for classification} 

Recently, \cite{zhao18} proposed the following definition of a breakdown point that is suitable for classification, calling it ''angular breakdown point'' since it is based on the angle between the decision hyperplane of the original coefficient and the one estimated on a contaminated sample. The following definition stems from \cite[Def. 1]{zhao18} and assumes linear classifiers. 

\begin{Def}[Angular breakdown point for classification] The \textbf{(population) angular breakdown point for classification} is given by \begin{equation} \label{angbdp} \epsilon(\beta,Z_n)=\min\left\{ \frac{m}{n} \ \bigg| \ \hat \beta(Z_n^m) \in S^-\right\}, \ \ \ S^-=\{\tilde \beta \ | \ \tilde \beta^T \beta \le 0\}. \end{equation} \end{Def}

\cite[Def. 1']{zhao18} also proposed a sample counterpart of this breakdown point where $\beta$ is replaced by $\hat \beta(Z_n)$ and therefore $S^-$ by $\hat S^-$ with the respective replacement. The angular breakdown point indicates that modifying more than $\epsilon(\beta,Z_n)$ of the sample $Z_n$ by arbitrary points can induce an angle between the original decision hyperplane and the hyperplane of the coefficient corresponding to the contaminated sample of at least $\pi/2$, leading to very low discriminative power if the classifier corresponding to $\hat \beta(Z_n)$ was sufficiently well. This setting has been extended to multi-class classification in \cite{qian}.

\section{Outliers and breakdown for ranking with linear scoring functions} \label{motivsec} 

As a motivation, we consider continuous ranking problems where the responses are continuously-valued (taking values in wlog. the whole space $\R$) in this section. 

\subsection{Why neither the regression nor the classification breakdown point work} 

We start by proving a counterpart of \cite[Prop. 3.1]{zhao18} showing that the finite-sample breakdown point in Eq. \ref{fsbdp} in also not reasonable in the ranking context. To this end, let the objective function of regularized continuous ranking with linear scoring functions be given by \begin{center} $\displaystyle L_{\lambda,n}(b,\beta,Z_n)=\frac{1}{n(n-1)}\mathop{\sum \sum}_{i<j} L((Y_i-Y_j)(s_{b,\beta}(X_i)-s_{b,\beta}(X_j)))+J_{\lambda}(\beta) $ \end{center} where $s_{b,\beta}(x):=x\beta+b$ so that $s(X_i)=:\hat Y_i$ is a parametric scoring function for some optional intercept $b$ with $|b|<\infty$, a loss function $L$ as introduced in Sec. \ref{prelim} and a regularizer $J_{\lambda}(\beta)$ satisfying \begin{equation} \label{penalty} \begin{gathered} \textbf{i)} \ J_{\lambda} \ge 0, J_0 \equiv 0, \ \ \ \textbf{ii)} \ J(\beta)=0 \Longleftrightarrow \beta=0_p, \\ \textbf{iii)} \ J(-\beta)=J(\beta), \ \ \ \textbf{iv)} \ J(\beta) \overset{||\beta|| \rightarrow \infty}{\longrightarrow} \infty, \end{gathered} \end{equation} where $0_p$ is the vector of length $p$ containing only zeroes. The fourth property is also known as coercivity (e.g. \cite{werner06}). The regularizer encourages sparse models and therefore does not take the intercept $b$ into account. 

Then, having a contaminated sample $Z_n^m$, we can rewrite the objective as \begin{equation} \label{objregrank} \begin{gathered} L_{\lambda,n}(b,\beta,Z_n^m):=\left[\frac{1}{n(n-1)}\mathop{\sum \sum}_{i \ne j, i,j \in I} L((Y_i-Y_j)(s_{b,\beta}(X_i)-s_{b,\beta}(X_j)))+J_{\lambda}(\beta)\right] \\ +\frac{1}{n(n-1)}\mathop{\sum \sum}_{i<j, i,j \in I^0} L((Y_i^0-Y_j^0)(s_{b,\beta}(X_i)-s_{b,\beta}(X_j))) \\ +\frac{1}{n(n-1)}\mathop{\sum \sum}_{i<j, i \in I^0, j \in I} L((Y_i^0-Y_j)(s_{b,\beta}(X_i^0)-s_{b,\beta}(X_j))) \\ =:G_{\lambda,n}(\tilde \beta,Z_{n-m})+F_n(\tilde \beta,Z_m^0)+H_n(\tilde \beta,Z_{n-m},Z_m^0) \end{gathered} \end{equation} for $\tilde \beta=(b,\beta)$ and where $Z_m^0$ denotes the contaminated part of the sample, i.e., $Z_m^0=\{(X_i^0,Y_i^0), i=1,...,m\}$, $Z_{n-m}$ is the clean part of the sample and $I^0$ denotes the indices of the contaminated instances w.r.t. $Z_n$, i.e., $I^0$ is an $m-$subset of $\{1,...,n\}$ such that $I \cup I^0=\{1,...,n\}$ and $I \cap I^0=\emptyset$ for the indices $I$ of the clean instances w.r.t. $Z_n$. Note that \cite{zhao18} do not need the $H_n-$term since there are no interactions between clean and contaminated instances. 

\begin{Def} A sample $Z_n$ is \textbf{linearly inrankable} if there exists no linear scoring function (linear in $\beta$) $s_{b,\beta}(x)=x\beta+b$ such that we can perfectly replicate the ranking of the responses in $Z_n$, otherwise we call the sample \textbf{linearly rankable}. \end{Def} 

Graphically, the most simple case linear rankability can be easily understood in the sense that the sign of the differences of the responses along each axis coincides for all instances, so finding a coefficient with the correct sign in each entry ensures a perfect ranking. However, linear inrankability \textbf{is generally not given} if this property is violated for some axes. 

\begin{ex} Consider the sample $((1,1),1)$, $((0,3),2)$, $((3,2),3)$. The responses clearly do neither strictly monotonically increase with increasing $X_{\cdot,1}$ nor with increasing $X_{\cdot,2}$. However, for $\beta:=(1,1)$ and arbitrary but finite $b$, we have $\hat Y_1=2+b$, $\hat Y_2=3+b$, $\hat Y_3=5+b$, so the ranking is perfect. \end{ex}

This example points out that linear inrankability does not only depend on some strict monotonicity of the responses w.r.t. some variable but also on the variables themselves (unless the strict monotonicity is satisfied along all axes). This is a first motivating aspect which makes the angular breakdown point from \cite{zhao18} inappropriate for ranking. Let us state the following counterpart to \cite[Prop. 3.1]{zhao18}. 

\begin{lem} \label{infcoef} Let $L$ be a nonnegative loss function with $L(0)<\infty$ and let the assumptions in Eq. (\ref{penalty}) be true. \\
\textbf{a)} For $\lambda>0$, it holds that $||\hat \beta(Z_n)||<\infty$ and $||\hat \beta(Z_n^m)||<\infty$ for any $Z_n$ and $Z_n^m$. \\
\textbf{b)} For $\lambda=0$, norm finiteness of the estimated coefficient cannot be guaranteed.
\end{lem}

\begin{bew}\textbf{a)} Since $J_{\lambda}(\beta) \rightarrow \infty$ as $||\beta|| \rightarrow \infty$, we just have to show that there exists a $\beta$ with $||\beta||<\infty$ such that the objective is finite. This is true since $\beta=0_p$ leads to a finite loss and $J_{\lambda}(0_p)=0$, so there definitely exists an optimizer of $L_{\lambda,n}$ with finite norm, disregarding if we have $Z_n$ or $Z_n^m$. \\
\textbf{b)} Recall that ranking loss functions are based on the product of the differences of the responses resp. the fitted responses. Clearly, having infinite values, we face the problem that it is impossible to reasonably define something like ''$\infty-\infty$'' which would arise if for example $\beta=(\infty,-\infty)$ and $X_i=(1,1)$. However, the indicator loss function in the hard ranking loss can be simply rewritten as \begin{center} $ \displaystyle I(\{\{Y_i>Y_j\} \cap \{\hat Y_i<\hat Y_j\}\} \cup \{\{Y_i<Y_j\} \cap \{\hat Y_i>\hat Y_j\}\}), $ \end{center} so the loss will always be computable. If $||\beta||=\infty$ but if $\hat Y_i=s_{b,\beta}(X_i)$ and $\hat Y_j=s_{b,\beta}(X_j)$ are computable, we can indeed get an infinite norm solution. For example, consider univariate predictors so that $\beta=\infty$ leads to $s_{b,\beta}(x)=\pm \infty$ for $\sign(x)=\pm 1$ (let $\sign(0):=0$). Then, if the signs coincide, which is true if for example all $X_i>0$, we just get the loss $L(0)$ as for random guessing which proves the statement since there is no evidence that there exists a coefficient with finite norm which can beat each coefficient with infinite norm.  \begin{flushright} $_\Box$ \end{flushright} \end{bew} 

\begin{rem} \label{randguess} The potential incomputability of the scores for infinite coefficient components is a severe problem. Therefore, we propose to treat the whole ranking model as random guessing (i.e., $\beta=0_p$) if there exists any $i$ such that $s_{b,\beta}(X_i)$ is not computable. \end{rem}

To illustrate the second part further, let us look at the following very simple example.

\begin{ex} \label{infperf} Consider one of the most simple cases that one could imagine, i.e., we have the sample $(-1,-1)$, $(1,1)$. Then the coefficient $\hat \beta=\infty$ will produce a perfect ranking which cannot be beaten by any other coefficient. \end{ex} 

Lemma \ref{infcoef} indicates that the usual finite-sample breakdown point in Eq. \ref{fsbdp} is insufficient for measuring the robustness of regularized ranking problems since any contamination keeps the norm of the estimated coefficient finite if the assumptions in the lemma hold. As for the angular BDP for classification introduced in \cite{zhao18}, we similarly can conclude that it is inappropriate for ranking if the variables are scaled differently or if they take values in different spaces. We provide a simple but illustrative example. 

\begin{ex} Let the sample $((5,0.2),0.9)$, $((6,0.3),1.2)$, $((1,0.1),0.3)$ be given an let $\beta=(0.1,2)$ be the true coefficient (wlog. let $b=0$). Then for $\tilde \beta=(0.2,-1)$ we have $\beta \tilde \beta<0$ but the ordering of the predictions w.r.t. $\tilde \beta$ is still correct. \end{ex}

\begin{rem}[Interpretation of the classical BDP] Let us additionally highlight the fact that one has to be very cautious in interpreting the classical BDP. Focusing on the values $\pm \infty$ can be highly misleading since it does not fully reflect the real meaning of the BDP. 

Let us consider the median whose BDP is given by $n/2$ for even $n$ if the two candidate observations are averaged and $(n+1)/2$ for uneven $n$. Let $n$ be uneven. Then the BDP does not indicate solely that we can get a median value of $\infty$ by manipulating $(n+1)/2$ of the observations but that we can find a manipulation of $(n+1)/2$ observations such that the median can take an arbitrary value. In other words, if we have access to $(n+1)/2$ observations, we have full control over the estimated median, and if one indeed would set the outliers to $\infty$, the median would also take the value $\infty$. 

Now, consider a linear regression estimator. Similarly, having control over one of the observations allows us to produce an arbitrary estimated regression coefficient, as the proof of \cite[Thm. 1]{alfons13} reveals. However, note that it is impossible to achieve an estimated coefficient with $||\beta||=\infty$ for arbitrary data. Just consider the case $p=1$. By letting for example the response of the mostright observation grow, the regression coefficient clearly would also grow, but in the limit case that the response takes the value $\infty$, the result would be that any coefficient leads to a loss of $\infty$ due to the original data points for which an infinite coefficient would predict an infinite value, therefore one cannot enforce a coefficient with infinite norm. \end{rem}

The scope of this work is to define a BDP for ranking such that a breakdown enables the worst-case ranking which reverts the true ordering. 

\begin{rem}[Non-linear scoring functions] The restriction to linear scoring functions (i.e., linear in $x$) is not necessary since sign-reverting all components of $\beta$ would clearly also revert the ordering for any scoring function of the form $s_{b,\beta}(x)=f(x)\beta+b $ where $f: \R^p \rightarrow \R^{p'}$ maps the regressors from the original regressor space $\mathcal{X}$ to some feature space $\mathcal{X}' \subset \R^{p'}$ where we allow $p \ne p'$ which refers for example to very natural situations like facing categorical regressors whose encoding enlarges the column number of the regressor matrix. The respective outlier configurations that we provide in the remainder then have to be concentrated on regions where the score is strictly monotonic w.r.t. the original coefficient. 

Note that the reduction to linear scoring functions is done for the sake of simplicity and illustrativeness and no restriction (as long as our scoring functions are still linear in the parameter $\beta$ which evidently is the case) since one could essentially approximate any non-linear scoring function by piece-wise linear scoring functions. The case of kernel-based scoring functions will be discussed in Sec. \ref{furthersec}. \end{rem} 

\subsection{The order-inversal breakdown point for ranking} 

Before we state the definition of our OIBDP for ranking, we proceed along the same lines as \cite{zhao18} and study the effect of outliers. Analogously, we first consider a single outlier, wlog. $I^0=\{1\}$, i.e., we have the data set $\tilde Z_n=Z_{n-1} \cup \{Z_1^0\}$ where $Z_1^0=(X_1^0,Y_1^0)$ is some contaminated instance. Let \begin{center} $ \displaystyle (\hat b, \hat \beta(\tilde Z_n))=\argmin_{(b,\beta): |b|<\infty}(L_{\lambda,n}(\tilde \beta,\tilde Z_n))$. \end{center} Let the loss function satisfy $\lim_{u \rightarrow -\infty}(L(u))=\infty$ for illustration. Then, for $||X_1^0|| \rightarrow \infty$ and $||Y_1^0|| \rightarrow \infty$, the loss diverges for each instance such that $(\hat Y_i-\hat Y_i)(Y_i-Y_j) \rightarrow -\infty$. Here, we have to ensure that $(s_{\hat b,\hat \beta}(X_i)-s_{\hat b,\hat \beta}(X_1^0))(Y_i-Y_1^0) \ge 0$ for all $i \ne 1$. These conditions correspond to the constraint set \begin{center} $ \displaystyle \check S_{Z_1^0}^+:=\bigcap \check S_{Z_1^0}^+(i), \ \ \ \check S_{Z_1^0}^+(i):=\{\beta \ | \ (s_{b,\beta}(X_i)-s_{b,\beta}(X_1^0))(Y_i-Y_1^0) \ge 0\} $ \end{center} which mimicks the set $S_{Z_1^0}^+$ from \cite[p. 8]{zhao18}, so that we can formulate the optimization problem as \cite[Eq. (3.3)]{zhao18} as \begin{center} $ \displaystyle \min_{|b|<\infty, \beta \in \check S_{Z_1^0}^+}(G_{\lambda,n}(\tilde \beta, Z_{n-1}))$ \end{center} with $G_{\lambda,n}$ from Eq. \ref{objregrank}. Clearly, the set $\check S_{Z_1^0}^+$ is a cone as an intersection of cones which is true since for $\beta \in \check S_{Z_1^0}^+(i)$ it holds that $c\beta \in \check S_{Z_1^0}^+(i)$ for any $c \ge 0$. For the case of $m$ outliers, we similarly get the constraint set \begin{center} $\displaystyle \check S_{Z_m^0}^+=\bigcap_{i \in I^0} \check S_{Z_i^0}^+. $ \end{center} 

In the ranking setting, we get interesting insights into the constraint set due to the fact that even in the single outlier case, \textbf{the pairwise nature of ranking problems lets the outlier act globally}, in contrast to the classification setting from \cite{zhao18} where each outlier acts locally. We now distinguish between $X-$ and $Y-$outliers. 

\textbf{i) $Y-$outliers:} If we have only outliers in the response, all $X_i$ are maintained but some $Y_i$ are contaminated such that one observes $Y_i^0$, i.e., the outlier set contains instances $(X_i,Y_i^0)$ for $i \in I^0$. Let wlog. $Y_1^0$ be the only outlier. Then, letting $||Y_1^0|| \rightarrow \infty$ can cause the data to be (linearly) inrankable. To see why this is only possible but not guaranteed, consider the most simple case that $p=1$. Let $Y_{i_1} \le ... \le Y_{i_n}$ for pair-wise different $i_j \in \{1,2,..,n\}$. Having linearly rankable data, i.e., $Y_j>Y_i \Leftrightarrow X_j>X_i$, replacing $Y_{i_1}$ with an extreme negative outlier, i.e., letting $Y_{i_1}^0 \rightarrow -\infty$, does not alter the linear rankability of the data, therefore, is not affecting the quality of the estimated coefficient so that its sign is maintained. The same holds for $Y_{i_n}^0 \rightarrow \infty$. On the other hand, for general $p$, if one of the intermediate responses is replaced by an extreme outlier, the set $\check S_{z_1^0}^+$ breaks down to $\{0_p\}$.
 \ \\

\textbf{ii) $X-$outliers:} Extreme $X-$outliers can have a similar effect. $X-$outliers solely affect the predictors, so the outliers are given by $(X_i^0,Y_i)$ for $i \in I^0$, i.e., the attacker cannot alter the responses. Consider again the example with $p=1$ as above and let $||X_{i_1}^0|| \rightarrow -\infty$. Again, this has no effect since any positive coefficient still produces a perfect ranking, similarly when replacing the regressor $X_{i_n}$ corresponding to the largest response $Y_{i_n}$ by an extreme positive outlier. 

We learned from this discussion that extreme outliers combined with an unbounded loss function are prone to induce zero coefficients $0_p$. As argued in the previous subsection, reverting the sign of some components of the coefficient does not guarantee any effect on the ranking quality unless all coefficient components are sign-reverted. Even this does not guarantee an inverted ordering but guarantees that the predicted ordering cannot be perfect anymore. Taking all these arguments into account, we now state the following definition for the OIBDP for ranking. 

\begin{Def}[order-inversal breakdown point for ranking] \label{bdprankdef}  \textbf{a)} The \textbf{population order-inversal breakdown point for ranking} is defined by \begin{center} $ \displaystyle \check \epsilon(\beta,Z_n):=\min\left\{\frac{m}{n} \ \bigg| \ \hat \beta(Z_n^m) \in S_{\cap}^-\right\}, \ \ \ S_{\cap}^-:=\bigcap_{j: \beta_j \ne 0}\{\tilde \beta_j \ | \ \tilde \beta_j \beta_j<0\} .$ \end{center}
\textbf{b)} The \textbf{sample order-inversal breakdown point for ranking} is defined by \begin{center} $ \displaystyle \check \epsilon(\hat \beta,Z_n):=\min\left\{\frac{m}{n} \ \bigg| \ \hat \beta(Z_n^m) \in \hat S_{\cap}^-  \right\}, \ \ \ \hat S_{\cap}^-:=\bigcap_{j: \hat \beta_j(Z_n) \ne 0}\{\tilde \beta_j \ | \ \tilde \beta_j \hat \beta_j(Z_n)<0\} .$ \end{center} \end{Def} 

\section{Asymptotic bounds for the breakdown point of the hard continuous ranking problem} \label{hardcontsec} 

In this section, we prove asymptotic upper bounds for the OIBDP (w.r.t. the loss functions) for the hard continuous ranking problem.

\subsection{Univariate case} 

We start by formulating the following assumption.

\begin{ass} \label{infass} Let $L$ be a continuous and strictly monotonically decreasing function with $\lim_{u \rightarrow \infty}(L(u))=0$ and $\lim_{u \rightarrow -\infty}(L(u))=\infty$. \end{ass}

These unbounded loss functions arise once convex surrogates are used, for example in RankBoost (\cite{freund}), RankingSVM (\cite{herb}, \cite{joachims02}) or the p-Norm-Push (\cite{rudinc}). 

\begin{rem}[Ties] \cite{zhao18} also respect the case of zero coefficients. %which leads to a lower bound loss that differs from $G_{\lambda,m}$. 
There are indeed situations where ties in ranking problems may get an individual loss, for example when computing the empirical AUC, see e.g. \cite{agar}, but for simplicity, we always assume that ties lead to a zero summand in the loss in this work. \end{rem}

\begin{rem}[Non-zero assumption] We always assume that the responses for the original data are mutually different in continuous ranking problems resp. that the responses in bipartite and $d-$partite ranking problems are not all the same so that the true coefficient is never $0_p$, so in other words, we exclude the issue of an \textbf{inlier breakdown} which is for example relevant for scale estimators, see e.g. \cite{huber}, \cite{hampel}. Together with the assumption of linear rankability with the original coefficient that we use in the proofs, this is the counterpart of assuming that the points are in general position.  \end{rem}

\begin{rem}[Immunization against particular regularization terms] We will not consider the regularization term directly in the following lemmas and theorems. Having $\lambda>0$ which defines a feasible set of the form $B_{q,c_{\lambda}}:=\{\beta \ | \ ||\beta||_r \le c_{\lambda}\}$ for some $0<c_{\lambda}<\infty$ and some $r>0$, we can always assume that we project the true coefficient onto this set by standardizing all components uniformly (which does not alter the ranking) since the issue of a sign-reversal does not depend on the magnitude of the respective coefficient components. We will discuss to the case $r=0$ in Sec. \ref{discusssec}. \end{rem}

\begin{lem} \label{infbdplemma} Let $p=1$. For the hard ranking problem with the loss function $L(u)=I(u<0)$, the sample and population OIBDP for ranking is \begin{equation} \label{bdpcond} \frac{\check m}{n}, \ \ \ \check m=\min\left\{m \ \bigg| \ m(n-m)+\frac{m(m-1)}{2}>\frac{(n-m)(n-m-1)}{2} \right\}  \end{equation} and asymptotically, the BDP is given by $1-\sqrt{0.5}$. \end{lem}

\begin{bew} The proof is given for the population version, the sample version is proven completely analogously. Assume wlog. that $\beta>0$. Consider the worst-case outliers shown in Fig. \ref{outliers}.  
\begin{figure}[H] \begin{center}  \includegraphics[width=4cm]{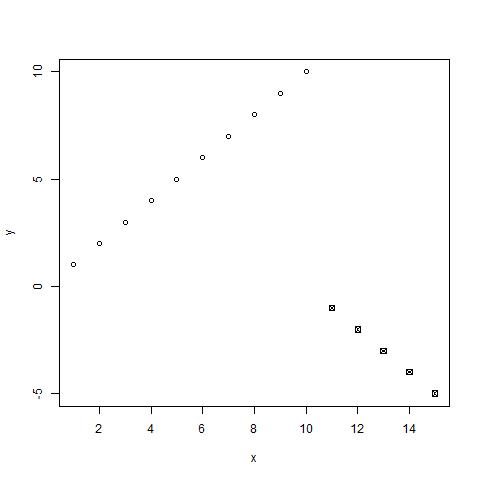}  \end{center} \caption[Worst-case outliers for $p=1$]{Worst-case outliers for $p=1$} \label{outliers} \end{figure} 

A correct ranking of two instances does not suffer a loss while any incorrect ranking suffers the same loss. This makes it impossible to achieve a breakdown by letting the response of one single outlier tend to $-\infty$ for all $n$ that are reasonably high ($\ge 4$). Observe that, due to symmetry, we have $n(n-1)/2$ effective pairwise comparisons and that a single outlier like the rightmost point in Fig. \ref{outliers} leads to $(n-1)$ misrankings for a coefficient $\beta>0$. Consider to add one outlier. Then, comparing each of the $(n-2)$ non-contaminated instances with one outlier leads to $(n-2)$ misrankings for $\beta>0$ but since the ordering of the outliers is also incorrect, we get a total of $2(n-2)+1$ misrankings. Now, let $m \ge 1$ outliers be given. Then, we get a total of $m(n-m)+m(m-1)/2$ misrankings for $\beta>0$ while the number of misrankings that we make for $\beta<0$ is evidently given by $(n-m)(n-m-1)/2$ since every pair of clean observations of which we have $(n-m)$ ones is misranked, so the number of outliers $\check m$ that we require for a breakdown is as stated in Eq. \ref{bdpcond}. 
\begin{figure}[H]  \begin{center}  \includegraphics[width=4cm]{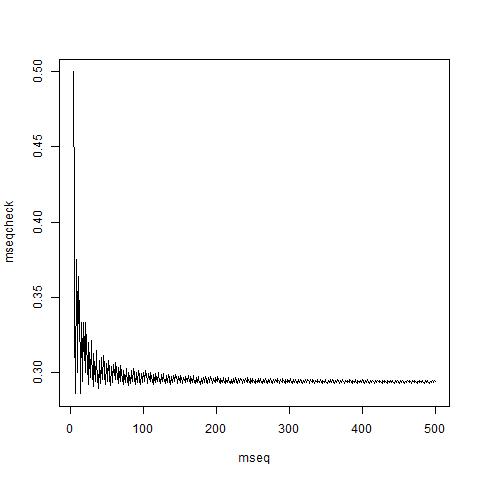} \caption[OIBDP for hard ranking for $p=1$]{OIBDP for hard ranking for $p=1$} \label{hardrankbdpplot}  \end{center} \end{figure}

Fig. \ref{hardrankbdpplot} shows the sample OIBDP for ranking for $n \in \{4,5,...,500\}$. Note that the maximal sample BDP is attained for $n=4$ where $\check m=2$ and the minimal sample BDP is attained for both $n=7$ and $n=14$, namely $\check m/n=2/7$.

As for the asymptotic setting, we set $m=cn$ and solve the inequality given in \ref{bdpcond} for $c$. \begin{center} $ \displaystyle cn^2-c^2n^2+\frac{c^2n^2}{2}-\frac{cn}{2} \overset{1}{>} \frac{(n-cn)(n-cn-1)}{2}=\frac{n^2-2cn^2+c^2n^2+cn-n}{2} $ \end{center} \begin{center} $ \displaystyle \overset{n>0}{\Longleftrightarrow} n\left[-c^2+2c-\frac{1}{2}\right]+[-c+0.5] \overset{!}{>} 0 . $ \end{center} Asymptotically, we just require that the value in the bracket of the left hand side is positive. We can easily conclude that this holds for $c>1-\sqrt{0.5}$, so this value is the sharp asymptotic upper bound for the OIBDP. \begin{flushright} $_\Box$ \end{flushright} \end{bew}

\begin{rem} Note that this result equals the asymptotic BDP of the Hodges-Lehmann estimator (see \cite[Sec. 11]{hodges}). This is not surprising since the Hodges-Lehmann estimator is given as the median of the set of all possible pairs of univariate samples. In order to achieve a breakdown of such an estimator, at least the half of the underlying observations which, in case of the Hodges-Lehmann estimator, are pairwise comparisons, have to be contaminated. This equals our setting since for the Hodges-Lehmann estimator, we need sufficiently many outliers in order to have control over more than the half of the pairs, i.e., we have to control more than $n(n-1)/4$ pairs. More precisely, the sum of the $m(m-1)/2$ outlier-outlier-pairs and the $m(n-m)$ outlier-non-outlier-pairs has to be greater than $n(n-1)/4$.% i.e.,  \begin{equation} \label{bdpcond2} \frac{\check m}{n}, \ \ \ \check m=min\left\{m \ \bigg| \ m(n-m)+\frac{m(m-1)}{2}>\frac{n(n-1)}{4}\right\}.  \end{equation} A similar computation as in the proof of Lemma \ref{infbdplemma} %shows that \begin{center} $ \displaystyle cn(1-c)n+\frac{cn(cn-1)}{2}=c(1-c)n^2+\frac{c^2n^2}{2}-\frac{cn}{2} \overset{!}{>} \frac{n^2}{4}-\frac{n}{4} $ \end{center} \begin{center} $ \displaystyle \Longleftrightarrow n^2\left[c(1-c)+\frac{c^2}{2}-\frac{1}{4}\right] \overset{!}{>} n\left[\frac{c}{2}-\frac{1}{2}\right] \overset{n>0}{\Longleftrightarrow} n\left[c-\frac{c^2}{2}-\frac{1}{4}\right] \overset{!}{>}\frac{c}{2}-\frac{1}{4}. $ \end{center} Asymptotically, we again have a again leads to an asymptotic BDP of $1-\sqrt{0.5}$. 
\end{rem} 

%\begin{rem} Other robust estimators that consider pairs of observations are the estimators $S_n$ and $Q_n$ from \cite{croux93}. \end{rem}

We do not consider unbounded loss functions here since we cover this case with Thm. \ref{supbdpinf} in the next subsection. Before we proceed with bounded loss functions, we argue why it suffices to consider indicator functions here. 

\begin{rem}[Reduction to indicator loss functions] Assume that the loss function is bounded, i.e., $\lim_{u \rightarrow -\infty}(L(u))=C_l<\infty$. Then, we can obviously generate losses which are (close to) $C_l$ by using the outlier scheme introduced above. However, since BDP computations have to consider all possible data configurations, we cannot exclude that the original data are so that for a broken coefficient, one suffers a loss $C_l-\epsilon$ for each comparison of original instances. More precisely, there is no guarantee that the loss suffered by a non-broken coefficient on the outliers is (considerably) greater than the loss suffered by a broken coefficient on the original data points. If the maximum $C_l$ is attained for a finite $u$, we always find original data such that we exactly suffer a loss of $C_l$ for a broken coefficient so the loss suffered on the outliers for the original coefficient equals the loss suffered on the original data points for the sign-reverted coefficient as for the indicator loss function. Therefore, we can restrict ourselves to the case of the indicator loss function where only the number of correct rankings resp. misrankings is taken into account. \end{rem} 

\subsection{Multivariate case} 

We considered the case $p=1$ separately since the arguments and results for $p>1$ are different for ranking. Now, the question arises if a breakdown in the sense of the OIBDP for ranking can always be achieved by manipulating $m<n$ data points, disregarding the particular configuration of the original data and the dimension.

\subsubsection{Unbounded loss function}

\begin{thm}\label{supbdpinf} Let $L$ satisfy Assumption \ref{infass}. Then, the upper bound for the sample and population OIBDP for ranking is $(p+1)/n$ provided that $1<p<n-1$, $p/n$ for $p=n-1$, $1/n$ for $p=1$ and not existent otherwise. \end{thm}

\begin{bew} Let us first illustrate our outlier configuration for $p=2$ and as usual, we only prove the population variant. Let wlog. be $\beta_1, \beta_2>0$ and let the original data points be linearly rankable according to $\beta$. Consider $X'=(\max_i(X_{i1}),\max_i(X_{i2}))$ and let $X^{(1)}=(X_1'+1,X_2')$ and $X^{(2)}=(X_1',X_2'+1)$. Set $Y'<\min_i(Y_i)$ and let $Y^{(1)}=Y^{(2)}<Y'$. This special configuration ensures that along each axis, any coefficient $\beta$ with positive components will produce a misrankings w.r.t. the outliers and, in addition, that there is no consistency with the original data since it is guaranteed that the response for all outliers would be greater than the response for all original variables according to any $\beta$ with positive components. Therefore, letting $Y^{(1)}, Y^{(2)} \rightarrow -\infty$, we produce an unlimited loss unless $\beta_1, \beta_2<0$. If any component, say, the first component of the original coefficient is negative, use $X^{(1)}=(X_1'-1,X_2')$ and proceed along the same lines. 

This strategy obviously is applicable to the general case $p>2$, requiring at most $(p+1)$ outliers. In the special case $p=n-1$, we just need $p$ instead of $p+1$ outliers by using the last remaining original data point as starting point for the construction of the $p$ outliers. The special case $p=1$ does not require a starting point.

As for the case $p \ge n$, consider for simplicity again $p=2$ and let $\beta_1, \beta_2>0$ and let $X_{11}<X_{21}$, $X_{12}<X_{22}$ and $Y_1<Y_2$. Regardless of the outlier configuration, one can only enforce the sign-inversal of one component. Even if one modifies $(X_2,Y_2)$ by letting $Y_2 \rightarrow -\infty$, any coefficient with $\beta_1>0$ and $\beta_2<0$ resp. $\beta_1<0$ and $\beta_2>0$ produces a perfect ranking provided that the negative component dominates here. Enforcing a sign-reversal of both components stays impossible and carries over to higher dimensions $p>2$. \begin{flushright} $_\Box$ \end{flushright} \end{bew}

\begin{rem} \label{gradrem} Note that although we cannot enforce multiple components to be sign-reverted by a single outlier, for particular algorithms this may not hold. When computing the BDPs, we considered all variables separately by our outlier schemes. This guarantees that our results provide conservative but valid upper bounds for the BDP which are sharp for the situations assumed in the respective theorems and lemmas. However, from an algorithmic point of view, the true BDP may be considerably lower. This is true for example for gradient-based approaches which update multiple coefficient components at once by a joint gradient step so that a single outlier in a remote location may pull the estimated coefficient towards a broken coefficient. However, this is a property of the numerical procedure that intends to minimize the corresponding objective and no general property. One may alleviate this issue by for example considering single gradient steps like in Gradient Boosting (e.g., \cite{bu07}). \end{rem}

Thm. \ref{supbdpinf} has a very interesting consequence: In a high-dimensional setting where $p>n-1$, it is impossible to find outlier configurations that guarantee a breakdown of the estimator in the sense of the OIBDP for ranking! Even if the whole data set would be replaced by outliers, is can only be enforced that $n$ coefficients are sign-reversed. Also note that it is not unusual that the dimension enters the BDP which also appeared for example in the BDP of the Least Trimmed Squares (LTS) estimator introduced in \cite{rous84}, see also \cite{rous06} for its fast computation, given in \cite{rous05}, where however an increasing dimension leads to a decreasing BDP. The sparse variant SLTS (\cite{alfons13}) also has a dimension-independent BDP. 

Going back to our ranking setting, the asymptotic case has to take the behaviour of the predictor dimension into account. See Sec. \ref{discusssec} for further discussions.

\begin{cor} Asymptotically, we have to distinguish between four cases. \\
\textbf{i)} If $p$ is fixed, then the asymptotic breakdown point is zero. \\
\textbf{ii)} If $p=p(n)=b_nn$ such that $b_n \rightarrow b \in [0,1[$, the asymptotic breakdown point is $b$. \\
\textbf{iii)} If $p=p(n)=b_nn$ such that $b_n \rightarrow b \ge 1$, the asymptotic breakdown point does not exist, i.e., it is impossible to achieve a breakdown for ranking.\end{cor}

The third case appears in usual high-dimensional settings, see e.g. \cite{bu}. 

\subsubsection{Bounded loss function}

\begin{thm} \label{supbdpind} Let wlog. $L$ be the indicator loss function used in the hard ranking loss and let $p \ge 2$. Then, the upper bound for the OIBDP for ranking is given by \begin{equation} \label{kstarhard} \frac{m^*}{n}, \ \ \ m^*=1+pk^*, \ \ \ k^*=\min\left\{k \ \bigg| \ \frac{k(k+1)}{2}>\frac{(n-pk-1)(n-pk-2)}{2}\right\}.  \end{equation} This quantity always exists for $p \le n-1$. \end{thm}

\begin{bew} Let us again illustrate our idea for $p=2$. The problem is that when having a starting point $X'$ as in the proof of Thm. \ref{supbdpinf}, generating one outlier by altering one component may not suffice to ensure a breakdown if the original data points still dominate. Moreover, we have to guarantee that all components of the coefficient are sign-reversed. We propose the following outlier algorithm: 

\begin{algorithm}
$k=1$\;
\While{\text{No breakdown} $\wedge \ 1+p(k+1)<n$}{
\For{$j=1,...,p$}{
Generate an outlier around $X'$ on the $j-$th axis as in the proof of Thm. \ref{supbdpinf}\;
\If{\text{Breakdown}}{Stop}}
$k=k+1$\;
}
$m=1+pk$
\end{algorithm}

For illustration, let $p=2$ and $k=2$. Then we generate a further outlier on each axis by proceeding on the respective axis, i.e., if $X^{(1)}=(X_1'+1,X_2')$, the next outlier is $X^{(3)}=(X_1'+2,X_2')$ and $Y^{(3)}=Y^{(4)}<Y^{(1)}=Y^{(2)}$ for $X^{(2)}=(X_1',X_2'+1)$ and $X^{(4)}=(X_1',X_2'+2)$. Applying this strategy, we get $k(k+1)/2$ comparisons along each axis, so keeping the original sign of the corresponding coefficient component leads to $k(k+1)/2$ misrankings. In contrast, we still have $(n-2k-1)$ original data points which, in the worst case, cause $(n-2k-1)(n-2k-2)/2$ misrankings provided that the coefficient has at least one component with the original sign. 

Note that comparisons of original and contaminated data points are not informative. Let us elaborate this argument a bit further. By construction, the responses of the outliers are lower than the responses of the original data which makes their ranking prediction perfect if all components of the coefficient are sign-reversed. In this case, the loss suffered due to these $(n-2k-1)(2k+1)$ comparisons is zero, so such a coefficient indeed leads to $k(k-1)/2$ misrankings. On the other hand, the original coefficient induces $(n-2k-1)(2k+1)$ misrankings, but any other coefficient in between these two extreme cases potentially predicts the respective orderings perfectly, so we have to be conservative and assume this ''least favorable case'' (from the view of the attacker) that such a coefficient also achieves a loss of zero like the completely sign-reverted coefficient when comparing outliers and original data points. Therefore, a breakdown is guaranteed once $k$ is large enough such that \begin{center} $ \displaystyle \frac{k(k+1)}{2}>\frac{(n-2k-1)(n-2k-2)}{2} $ \end{center} which leads to the stated formula \ref{kstarhard} for $p=2$. 

In the general case $p>2$, we consider at most $p-$chunks of $k$ new outliers, i.e., $m=1+pk$, and by the same arguments, a breakdown occurs if \begin{center} $ \displaystyle \frac{k(k+1)}{2}>\frac{(n-pk-1)(n-pk-2)}{2}. $ \end{center} 

Clearly, there exist cases where such a $k^*$ does not exist. Here, we have to distinguish between two cases: \textbf{i)} $p \ge n$; \textbf{ii)} $p \le n-1$.  \\
\textbf{i)} This case has already been discussed in Thm. \ref{supbdpinf} where we concluded that it is impossible to guarantee a breakdown in such high dimensions. This evidently also holds for the case of bounded loss functions. \\
\textbf{ii)} A breakdown may be achieved before a $p-$chunk is complete. In the worst case, we can stop once $m=n-1$ since then, using the last remaining point as starting point, we can generate at least one outlier along each axis. In general, provided that $k^*$ exists, the true upper bound BDP therefore lies in the set \begin{center} $ \displaystyle \left\{\frac{1+p(k^*-1)+1}{n},..., \frac{1+pk^*}{n}\right\} $. \end{center} Note that there exist configurations in which $m^*=1+pk^*$ is indeed sharp which is true for example for $p=n-1$ as already discussed. \begin{flushright} $_\Box$ \end{flushright} \end{bew}

\begin{ex} To illustrate the case ii) in the proof above, consider the case $p=2$ and $n=8$. Generating an outlying starting point and two outliers along each axis leads to $m=5$ and $k=2$, but we have only three comparisons of outliers along each axis and three comparisons of original instances. The loss for each coefficient $\beta$ with $\beta_1, \beta_2>0$ is obviously greater than the loss for each coefficient with $\beta_1,\beta_2<0$, but there is no guarantee that such a sign-reversed coefficient would achieve a lower loss than a coefficient with only one sign-reverted component. However, adding one additional outlier according to our outlier scheme, disregarding on which of the two axes, leads to a breakdown since the number of comparisons between original data boils down to one, leading finally to $m^*=6$ instead of $m^*=7$. \end{ex}

\begin{cor} \label{supbdpindcor} The asymptotic upper bound for the OIBDP for ranking \\
\textbf{i)} is given by $p/(p+1)$ for fixed $p$, \\
\textbf{ii)} is given by 1 for $p=p(n)=b_nn$ with $b_n \rightarrow b \in ]0,1[$, \\
\textbf{iii)} does not exist for $p=p(n)=b_nn$ with $b_n \rightarrow b \ge 1$. \end{cor}

\begin{rem} \label{outlierrem} We do not exclude that there may exist even more sophisticated outlier schemes than ours which leads to a faster breakdown. However, our outlier scheme guarantees a breakdown, provided than $p$ resp. $p(n)$ is small enough, which would be very hard to show for outlier schemes than are not axis-based. An intuitive alternative that however does not work is given in the Appendix in Ex. \ref{outlierdice}. \end{rem}

\subsection{Expected OIBDP} 

Evidently, there are always pathological configurations of the original data that even immediately cause a breakdown or that hinder a breakdown but being extremely artificial, see Ex. \ref{compex}. %Consider for example the simple empirical variance estimator and let all observations be identical. Then, the estimated variance hits the boundary value zero, so a breakdown has occurred without even having placed one single outlier (or inlier). However, on the other hand, always assuming a ''worst-case configuration'' of the original data points from the view of the attacker, i.e., a configuration that makes it most challenging to find appropriate outliers to cause a breakdown, can also be misleading. We provide an example in the supplementary file. This is the starting point for a discussion on the actual meaning of a BDP and the question what a suitable reference data configuration to define the BDP can be. \ \\
Note that a comparable situation has already been investigated in \cite{horbenko12} who consider the expectation of the BDP w.r.t. the ideal distribution (which the original instances are assumed to follow), leading to a so-called expected BDP. Their motivation was to account for the fact that unfavorable configurations of the original data points only appear with low probabilities which helped them to get nonzero expected BDPs in the context of heavy-tailed distributions or when only partial equivariance is valid. 

However, in our setting, we have to be very cautious how to define an \textbf{expected OIBDP}. Evidently, assuming iid. instances $(X_i,Y_i)$ and computing the expectation w.r.t. the joint distribution would make no sense since iid. instances are all ranked equally in expectation. We indeed require a fixed design of the regressor matrix which, for every fixed $n$, assumes that observations $X_{n,i}, i=1,...,i_n,$ are given. Then, the responses are computed by $Y_{n,i}=X_{n,i}\beta+\epsilon_{n,i}$ for $\epsilon_{n,i} \sim F_{\epsilon}$ iid. for some centered distribution $F_{\epsilon}$. Therefore, the points $(X_{n,i},X_{n,i}\beta)$ are trivially linearly rankable but the points $(X_{n,i},Y_{n,i})$ do not necessarily be linearly rankable since this property depends on the realizations of the error terms. In the proofs, we always consider linear rankability w.r.t. the original coefficient $\beta$ which can be interpreted as taking the expectation of the data w.r.t. $F_{\epsilon}$. This motivates the following definition which mimicks \cite[Def. 3.2]{horbenko12}. 

\begin{Def}[Expected OIBDP for ranking] \label{expbdp}  Let $Z_n(\epsilon)$ be the sample consisting of the data points $(X_{n,1},Y_{n,1}(\epsilon_{n,1})),...,(X_{n,i_n},Y_{n,i_n}(\epsilon_{n,i_n}))$. \\
\textbf{a)} The \textbf{expected population order-inversal breakdown point for ranking} is defined by $\E_{\epsilon}[\check \epsilon(\beta,Z_n(\epsilon))].$
\textbf{b)} The \textbf{expected sample order-inversal breakdown point for ranking} is defined by $\E_{\epsilon}[\check \epsilon(\hat \beta,Z_n(\epsilon))].$\end{Def} 

\begin{rem} One has to be very cautious when considering the sample OIBDP (or general sample BDPs). This fact has been respected by \cite[Thm. 4]{zhao18} who indeed assume that the estimator does not yet break down on the original sample. As for ranking, our theoretical results on BDP bounds are founded on the expectation w.r.t. the error term, making the data linearly rankable w.rt. $\beta$. 

However, as we already highlighted, these linearly rankable original data points support the original coefficient most, i.e., any tie or other inconsistency reduces the required amount of outliers to let the estimator break down. This can indeed be problematic if the sample BDP is considered and if, maybe due to a large error variance, the estimated coefficient is insufficiently supported by the data. Let $p=1$ and let the original coefficient have a small magnitude. Then, a large error variance may cause the data points to oscillate with growing regressor value, so just imposing one outlier may already change the sign of the estimated coefficient. We think that such issues are prone for a low signal to noise ratio (SNR) and that there may exist something like a ''noise gap'' between the population and sample BDP variants. \end{rem}

\section{Hard binary and hard $d-$partite ranking problems} \label{hardbinsec} 

The goal of binary hard ranking problems is to find the correct ordering of all instances w.r.t. the probability to belong to class 1 for binary responses. In fact, one computes a real-valued scoring function so that the ordering of the scores is equivalent to an ordering of the respective probabilities. In $d-$partite ranking problems, one proceeds as in ordered logit regression by binning the scores. %, so one may think of a hard continuous ranking problem on the real-valued image space of the scoring function. Similarly, in $d-$partite ranking problems, one has categorical responses with $d$ ordered categories. As in ordered logit models, one has to compute a real-valued scoring function and $(d-1)$ suitable thresholds such that an instance is classified into class $k$ if the score is between the $(k-1)-$th and $k-$th threshold, using $-\infty$ as zero-th threshold and $\infty$ as $d-$th threshold. 
However, while an ordered classification model would be perfect if all instances get a score that is contained in the correct interval, hard $d-$partite ranking problems require that the ordering of all scores, and therefore also in the respective chunks, is correct. 

As for the OIBDP computation, let us distinguish between two cases: \textbf{i)} We have access to the real-valued pseudo-responses, so we are again in the usual continuous setting, making the results from Sec. \ref{hardcontsec} applicable; \textbf{ii)} The more realistic case is that we indeed only observe the categorical responses and that we only can produce outliers with responses in the respective discrete set. The main difference to the continuous case is that the outlier configuration becomes far less flexible. Let the loss function operate on the score scale, i.e., we use $s_{b,\beta}(X_i)$ instead of $\hat Y_i$ where the latter would be $\pm 1$ for binary ranking, otherwise we were in a classification setting.

\subsection{Unbounded loss}

\begin{cor} \label{bininflosscor} If the loss function satisfies Assumption \ref{infass}, the upper bound of the sample and population OIBDP for ranking is $(p+1)/n$ for $p \le n-2$, $p/n$ for $p=n-1$, $1/n$ for $p=1$ and not existent otherwise. \end{cor}

\begin{bew} Follows the same steps as in the proof of Thm. \ref{supbdpind} with the difference that we cannot explicitly produce $Y-$outliers on the score scale (since it is unobservable) but only on the response scale. For illustration, set $p=2$. We cannot produce extreme $Y-$outliers but we indeed can produce extreme $X-$outliers, so let wlog. $\beta_1, \beta_2>0$. We again use a starting point $X'=(\max_i(X_{i1}),\max_i(X_{i2}))$ and set $Y'=1$. Then, let $Y^{(1)}, Y^{(2)}=-1$ and $X^{(1)}=(X_1'+c_1,X_2')$, $X^{(2)}=(X_1',X_2'+c_2)$. Letting $c_1, c_2 \rightarrow \infty$ will induce an unbounded loss when comparing $(X',Y')$ with each of the both outliers, so $\beta_1, \beta_2<0$ is enforced. Clearly, this will also produce losses when comparing the original data pairs and the pairs with one original data point and one outlier, but they are finite. 

This strategy clearly also holds for $p=1$ and $p>2$ where a complete sign-reversal is again no longer guaranteeable for $p \ge n$. For $p=n-1$, we can use the last remaining original data point as starting point for the construction if its response is negative and an analogous construction otherwise. \begin{flushright} $_\Box$ \end{flushright} \end{bew}

\begin{rem} \label{bipkpartrem} Due to the ordering of the classes in $d-$partite ranking problems, a similar approach can be executed to show that the upper bound for the OIBDP for ranking in this setting is identical. One just has to set $Y'=d$ and set $Y^{(1)}=Y^{(2)}=1$. Letting the respective predictor components diverge, the breakdown should be achievable for any setting where $Y'$ has to be a higher label than $Y^{(1)}$ and $Y^{(2)}$, but however, using the extreme classes is the most logical configuration. \end{rem} 

\subsection{Bounded loss}

Let us now translate Lemma \ref{infbdplemma} and Thm. \ref{supbdpind} to the case of hard binary ranking with the indicator loss function. As already elaborated in the proof of Cor. \ref{bininflosscor}, we do not have access to the true underlying real-valued scores but only to the binary reponses which severely restricts the possible outlier configurations. Then, the idea is essentially the same as in AUC maximizing approaches (for example done in \cite{ataman} and \cite{rak04} for SVM-type and in \cite{clem08c}, \cite{clem10b} for tree-type approaches for ranking), i.e., the score for each instance of class 1 has to be higher than the score for each instance of class -1. 

\begin{lem} \label{infbdplemmabinary}  Let $p=1$. For the hard binary and $d-$partite ranking problem with the loss function $L(u)=I(u<0)$, the sample and population OIBDP for ranking is given by \begin{equation} \label{bdpcondbinary} \frac{\check m}{n}, \ \ \ \check m=2 \check k, \ \ \ \check k=\min\left\{k \ \bigg| \ k\lfloor \frac{n}{2}\rfloor+k\left(\lceil \frac{n}{2} \rceil-k\right)>\left(\lceil \frac{n}{2} \rceil-k\right)\left(\lfloor \frac{n}{2} \rfloor-k\right) \right\}  \end{equation} and asymptotically, the BDP is given by $1-\sqrt{0.5}$. \end{lem}

\begin{bew} Due to the discrete observable response space, we cannot apply the outlier scheme the we suggested in Fig. \ref{outliers}. We first need to identify the configuration of the original responses that supports the true coefficient most. Therefore, we wlog. assume that $\beta>0$ is the true coefficient. Then, for even $n$, the worst-case original data configuration (from the view of the attacker) is composed by a set of $n/2$ instances with $X_i>0$ and $Y_i=1$ and a set of $n/2$ instances with $X_j<0$ and $Y_j=-1$. For uneven $n$, either ''half'' contains $\lceil n/2 \rceil$ resp. $\lfloor n/2 \rfloor$ instances. Since we do not have a classification problem where one usually classifies all instances with a score greater than zero as class 1 instance and vice versa, we do not necessarily have to consider $X=0$ as ''boundary'', we do it just for the sake of easiness. 

\begin{figure}[H] \begin{center}  \includegraphics[width=4cm]{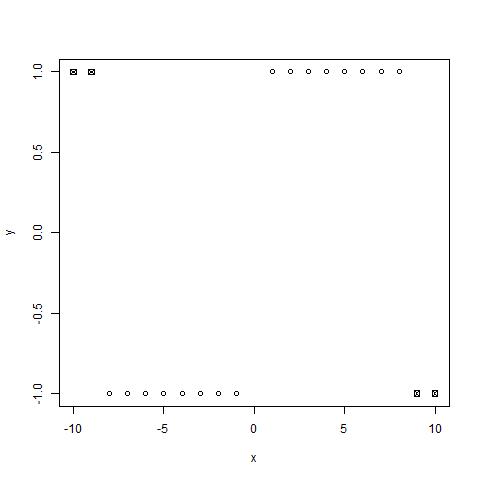}  \end{center} \caption[Worst-case outliers for $p=1$]{Worst-case outliers for $p=1$} \label{outliersbin} \end{figure}

Now, starting by replacing either the $m<n/2$ rightmost resp. leftmost instances (w.r.t. $X$) by outliers that keep the regressor value but switch the response sign would not lead to a breakdown if one only has access to the binary outcome, see Fig. \ref{outliersbin}. This is easily seen assuming $n (mod \ 4)=0$. If the $n/4$ rightmost instances are replaced as suggested, each $\beta>0$ will produce exactly $(n/4)^2$ misrankings which arise from comparing each of the $n/4$ original instances on the right half with the $n/4$ outliers. In contrast, any $\beta<0$ will produce $n^2/8$ misrankings by comparing each of the leftmost $n/2$ instances with the $n/4$ remaining instances on the right half. Note that one cannot compare the instances on the left half with the outliers since they all have the same response value. Replacing more than $n/4$ instances will even supply the original coefficient more whence we considered the case that around the half of the instances belong to either class as the worst case. 

Therefore, a reasonable outlier scheme is to start by replacing the leftmost and the rightmost instance simultaneously by outliers as illustrated in Fig. \ref{outliersbin}, i.e., by switching the sign of the response. Note that due to the binary response and the indicator loss function, $Y-$outliers and $X-$outliers can be regarded as being equivalent, so it suffices only to use $Y-$outliers, keeping the original regressor values. This first step induces, for even $n$, exactly $n/2$ misrankings by comparing the leftmost outliers with every instance of class -1 and additionally $(n/2-1)$ misrankings by comparing the remaining original instances on the right half with the rightmost outlier for $\beta>0$ (the sum will also be $(n-1)$ for uneven $n$). In contrast, the remaining $(n/2-1)$ original instances on the left half lead to $(n/2-1)(n/2-1)$ misrankings by comparing them with the $(n/2-1)$ original instances on the right half for $\beta<0$ (resp. $(\lfloor n/2 \rfloor -1)(\lceil n/2 \rceil-1)$ ones). 

Assuming that for step $k$, one has $k$ outliers on each side, i.e., a total of $2k$ outliers, one gets the requirement for $k^*$ stated in Eq. \ref{bdpcondbinary}. Asymptotically, we assume $k=cn$ and easily conclude that $c^*=\frac{1}{2}-\sqrt{1/8}$, so the asymptotic BDP is \begin{center} $ \displaystyle \frac{m^*}{n}=2c^*=1-\sqrt{\frac{1}{2}}. $ \end{center} 

We already argued in Rem. \ref{bipkpartrem} that outlier strategies for binary ranking problems are also applicable to $d-$partite ranking problems. Although one would have more flexible outlier schemes if the instance labels are diverse enough, for example, by considering ascending classes on the right half and descending classes on the left half which enables to produce more misrankings for the original coefficient by taking outlier-outlier-pairs into account, we would be essentially be in the same setting as in the binary ranking problem if one considers the ''least favorable'' configuration of the original data where the instances on the left half belong to class 1 and the ones on the right half to class $d$, making the upper bound of the OIBDP for the binary ranking problem also a sharp bound for the $d-$partite ranking problem. \begin{flushright} $_\Box$ \end{flushright} \end{bew}

Now, we consider the general case $p>2$. The proof of the following theorem and corollary can be found in the Appendix.

\begin{thm} \label{supbdpindbinary} Let $L$ be the indicator loss function and let $p \ge 2$. Then, the upper bound for the OIBDP for ranking for hard bipartite and hard $d-$partite ranking problems is given by \begin{equation} \label{kstarhardbinary} \frac{m^*}{n}, \ \ \ m^*=1+2pk^*, \ \ \ k^*=\min\left\{k \ \bigg| \ k(k+1)>\frac{(n-2pk-1)(n-2pk-2)}{2}\right\}.  \end{equation} This quantity always exists for $p \le n-1$. \end{thm}

\begin{cor} \label{supbdpindcorbin} The asymptotic upper bound for the OIBDP for bipartite and $d-$partite ranking \\
\textbf{i)} is given by $(2p^2-\sqrt{2}p)/(2p^2+1)$ for fixed $p$, \\
\textbf{ii)} is given by 1 for $p=p(n)$ with $p(n)/n \rightarrow b \in ]0,1[$, \\
\textbf{iii)} does not exist for $p=p(n)$ with $p(n)/n \rightarrow b \ge 1$. \end{cor}

\section{Localized ranking problems} \label{locsec} 

Localized ranking problems follow two goals, i.e., identifying the top $K$ instances and retrieve the ordering of the true or fitted top $K$ instances correctly. As for robustness analysis, we have to consider the OIBDP for ranking instead of the angular BDP for classification from \cite{zhao18} since the former one is stricter, so letting the ranking break down directly guarantees a breakdown of the classification due to the fixed number $K$ of class 1 instances (for $K<n/2$; the other case will also be discussed below). 

\subsection{Unbounded loss}

The proof can be found in the Appendix.

\begin{cor} \label{locinflosscor} If the loss function used for the ranking part satisfies Assumption \ref{infass}, the upper bound of the sample and population OIBDP for localized ranking is $(p+1)/n$ for $p \le K-2$, $p/n$ for $p=K-1$ and not existent otherwise. If the classification loss function satisfies Assumption \ref{infass}, the BDP is $(p+1)/n$ for $p \le K-1$, $p/n$ for $p=K$ and not existent otherwise. In either case, it is $1/n$ for $p=1$. \end{cor}

\subsection{Bounded loss}

As for the case of bounded loss functions, wlog. the indicator loss functions, we have to distinguish between a couple of cases, i.e., if $K \le n/2$ or $K>n/2$ and if the ranking part of the localized loss is based on $Best_K$ or on $\widehat{Best_K}$. \ \\

In this section, we require that the ranking of the true best $K$ instances is predicted correctly (see the discussion below Eq. \ref{locrankrisk}). The case of localizing on $\widehat{Best_K}$ can be found in the Appendix. 

\begin{lem} \label{infbdplemmaloc}  Let $p=1$. For the localized continuous ranking problem with 0/1-loss for classification and the indicator loss function for ranking where the latter is based on the true best instances with indices in $Best_K$, the sample and population OIBDP for ranking \\
\textbf{i)} is given by \begin{equation} \label{bdpcondlocbest} \begin{gathered} \frac{\check m}{n}, \ \ \ \check m=\min(K,\min\{k \ \bigg| \ \frac{n-K}{n} \cdot \frac{2(K-m)}{n}+\frac{(K-m)(K-m-1)}{2n(n-1)} \\ <\frac{n-K}{n} \cdot \frac{2m}{n}+\frac{1}{n(n-1)}\left[\frac{m(m-1)}{2}+m(K-m)\right] \}) \end{gathered}  \end{equation} for $K \le (n+m)/2$, \\
\textbf{ii)}  is given by \begin{equation} \label{bdpcondlocbest2} \begin{gathered}  \frac{\check m}{n}, \ \ \ \check m=\min\{k \ \bigg| \ \frac{n-K}{n} \cdot \frac{2n-K}{n}+\frac{(K-m)(K-m-1)}{2n(n-1)} \\ <\frac{n-K}{n} \cdot \frac{2m}{n}+\frac{1}{n(n-1)}\left[\frac{m(m-1)}{2}+m(K-m)\right] \}   \end{gathered} \end{equation} for $K \ge (n+m)/2$ and $K \le n-m$ provided that $n-m \ge (n+m)/2$, \\
\textbf{iii)} is given by Eq. \ref{bdpcond} in Lemma \ref{infbdplemma} where $n$ in the definition of $\check m$ is replaced by $K$ for $K \ge n-m$. \end{lem}

\begin{bew} \textbf{i)} Since we concentrate on the true top $K$ instances, we have to use the outlier scheme depicted in Fig. \ref{outliersloc} if $\beta>0$ is the true coefficient (for $\beta<0$, the responses of the $m$ rightmost instances would be moved upwards). \begin{figure}[H] \begin{center}  \includegraphics[width=4cm]{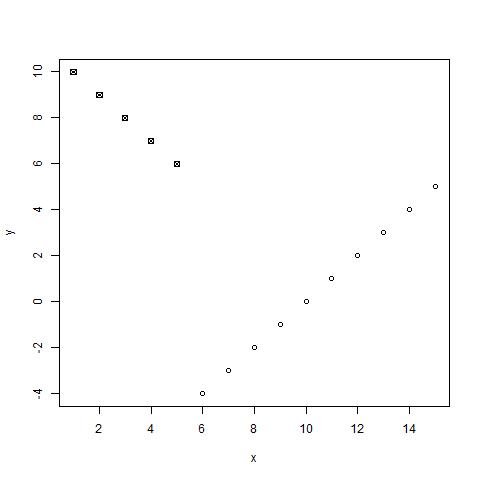}  \end{center} \caption[Worst-case outliers for the localized ranking problem for $p=1$]{Worst-case outliers for the localized ranking problem for $p=1$} \label{outliersloc} \end{figure} Then, any coefficient $\beta>0$ will misclassify $m$ instances, i.e., the outliers, while any coefficient $\beta<0$ will misclassify $(K-m)$ instances. As for the ranking part, any coefficient $\beta>0$ will produce both misrankings on the outliers as well as on every pair of an outlier an a clean top-$K-$instance, leading to $m(m-1)/2+m(K-m)$ misrankings, whereas any coefficient $\beta<0$ will produce misrankings on the non-outlier top-$K-$instances, i.e., $(K-m)(K-m-1)/2$ misrankings. Therefore, the statement in Eq. \ref{bdpcondlocbest} in the minimum-brackets is true. Moreover, if we had $m=K$ outliers, the loss for $\beta<0$ would be zero and since any positive coefficient produces a loss greater than zero, indicating that the BDP is not larger than $K/n$ which proves the outer minimum operator in Eq. \ref{bdpcondlocbest}. 

\textbf{ii)} The only difference in the case $K>n/2$ compared to the case $K \le n/2$ is that the broken coefficient does not necessarily cause $(K-m)$ misclassifications due to the fact that more than the half of the instances are labeled as class 1 instances. Consider the concrete example that $n=100$, $K=70$ and $m=10$, following the scheme in Fig. \ref{outliersloc}. Then $\beta<0$ does not lead to 60 misclassifications but just to 30 misclassifications, i.e., the rightmost $n-K=30$ instances. However, if we had $K=52$, we would not make $n-K=48$ misclassifications but just $K-m=42$ ones. Since $n-K<K-m$ iff $K>(n+m)/2$, Eq. \ref{bdpcondlocbest2} follows. 

\textbf{iii)} Once $m>n-K$, we cannot misclassify $m$ instances with the coefficient of the original sign but only $(n-K)$ ones (the leftmost $(n-K)$ ones in Fig. \ref{outliersloc}). The same is true for the broken coefficient where the rightmost $(n-K)$ instances in Fig. \ref{outliersloc} are misclassified. This lets the classification loss be equal for both cases and reduces the problem to the hard ranking problem on $Best_K$.  \begin{flushright} $_\Box$ \end{flushright} \end{bew}

\begin{cor} \label{infbdplemmaloccor} \textbf{i)} Asymptotically, a fixed $K$ would lead to a BDP of zero. For $K=n$, we get the asymptotic BDP of $1-\sqrt{0.5}$ as for hard ranking. \\
\textbf{ii)} For $K=K(n):=dn$ for $d \in ]0,d_0]$ with $d_0 \approx 0.6352578$, we can conclude that for $m=cn$, we have \begin{center} $ \displaystyle c^*=2-d-\sqrt{4-6d+5d^2/2} $ \end{center} which takes values in $]0,0.270514]$ and is strictly monotonically increasing w.r.t. $d$.  \\
\textbf{iii)} For $K=K(n):=dn$ for $d \in [d_0,d_1]$ with $d_1 \approx 0.773455$, we can conclude that for $m=cn$, we have \begin{center} $ \displaystyle c^*=1-\sqrt{-1+4d-5d^2/2}$ \end{center} which takes values in $[0.2265413,0.270514]$ and has its minimum at $d_1$.  \\
\textbf{iv)} For $K=K(n):=dn$ for $d \in [d_1,1]$, we have the asymptotic BDP $d(1-\sqrt{0.5})$.  \end{cor} 

\begin{bew} Statement i) is trivial. The formulae in statements ii) and iii) can be easily computed but we have to explain the value $d_0$. As we have seen in Eq. \ref{bdpcondlocbest2}, it depends on $m$ whether the asymptotic BDP corresponding to Eq. \ref{bdpcondlocbest} or to Eq. \ref{bdpcondlocbest2} applies, depending on $\min(K-m,n-K)$. Graphically, we search for the first intersection of both BDPs in dependence of $d$ (note that the second intersection is given at $d=1$). A numerical evaluation delivers the value $d_0$ above, see the black curve (asymptotic BDP from ii)) and the blue curve (asymptotic BDP from iii) which is only valid for at least $K>n/2$, therefore the growth for decreasing $d$ can be ignored) in Fig. \ref{asylocbdp}. Note that at $d_0$, the asymptotic BDP for both cases is given by around 0.270514 which equals $2(d_0-0.5)$ and can be explained by our argumentation in the proof of Lemma \ref{infbdplemmaloc} that the larger $m$ is in the case $K>n/2$, the lower is the number of misclassified instances for the broken coefficient, so according to the formula in Eq. \ref{bdpcondlocbest2}, we switch between both asymptotics once $n-dn>dn-cn$, i.e., once $c>2(d-0.5)$. For statement iv), we similarly have to find $d$ for which $dn>n-cn$ holds (i.e., $K>n-m$). Again, by numerical evaluation where we search for the intersection of the purple curve and the black line in Fig. \ref{asylocbdp} (for illustration, we also added the graph of $d \mapsto 1-d$ (red line) which intersects the blue curve at the same point) which is the case at $d_1$ where the BDP is exactly $1-d_1$, so increasing $d$ from here will cause the regime switch in the classification part of the right hand side, which results in the combined BDP curve in the right part of Fig. \ref{asylocbdp}. \begin{figure}[H] \begin{center}  \includegraphics[width=4cm]{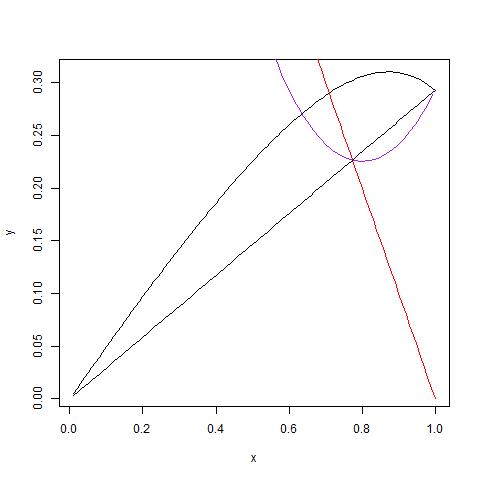}  \includegraphics[width=4cm]{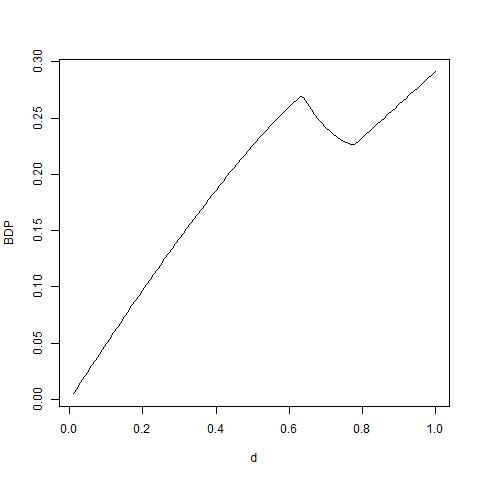}  \end{center} \caption[Asymptotic BDP for $K \le n/2$ and $K>n/2$]{Asymptotic BDP for $K \le n/2$ and $K>n/2$} \label{asylocbdp} \end{figure} \begin{flushright} $_\Box$ \end{flushright} \end{bew}

\begin{thm} \label{supbdpindloc} Let wlog. the indicator loss functions for both the classification and the ranking part be used and let $p \ge 2$. Then, the upper bound for the OIBDP for localized ranking is given by \begin{equation} \label{kstarloc}  \begin{gathered} \frac{m^*}{n}, \ \ \ m^*=\min(1+pk^*,K), \ \ \ k^*=\min\{k \ \bigg| \  \frac{n-K}{n}\frac{2\min(K-1-pk,n-K)}{n} \\ +\frac{(K-1-pk)(K-2-pk)}{2n(n-1)}<\frac{n-K}{n}\frac{2(1+pk)}{n}+\frac{k(k+1)}{2n(n-1)} \}   \end{gathered} \end{equation} for the case $K<n-m$. For $K>n-m$, we get the same $k^*$ as in Eq. \ref{kstarhard} in Thm. \ref{supbdpind}. This quantity always exists for $p \le K-1$. \end{thm}

\begin{bew} Along the same lines as the proofs of Thm. \ref{supbdpind} and Lemma \ref{infbdplemmaloc}.  \begin{flushright} $_\Box$ \end{flushright} \end{bew}

\begin{cor} \label{supbdpindloccor} For the localized ranking problem where the localized ranking loss is computed on $Best_K$ and where both the classification and the ranking loss are indicator functions, we asymptotically conclude that  \\
\textbf{i)} the BDP is zero for $K$ and $p$ being fixed,  \\
\textbf{ii)} the BDP for $c \le 0.5$ tends to \begin{center} $ \displaystyle pc^*, \ \ \ c^*=\frac{4p-3pd}{p^2-1}-\sqrt{\frac{16p^2-28p^2d+12p^2d^2+4d-3d^2}{(p^2-1)^2}} $ \end{center} and for $c>0.5$, it tends to the same quantity provided that $c \le 2(d-0.5)$, to \begin{center} $ \displaystyle pc^*, \ \ \ c^*=\frac{2p-pd}{p^2-1}-\sqrt{\frac{4p^2d-4p^2d^2-8d+5d^2+4}{(p^2-1)^2}} $ \end{center} for $c \in [2(d-0.5),1-d]$ and to $dp/(p+1)$ otherwise. For $p \rightarrow \infty$, the break-even point is given by $d_0 \approx 0.6923$ and the second break-even point tends to 1 for growing $p$.  \\
\textbf{iii)} the BDP tends to the asymptotic BDP for hard ranking, i.e., $p/(p+1)$, for $d=1$.  \\
\textbf{iv)} the BDP does not exist for $p=p(n)$ with $p(n)/n \rightarrow b \ge d$. \end{cor}

\section{Other ranking problems \label{othersec}} 

In this section, we briefly address the instance ranking problems that we did not consider so far. 

Weak ranking problems (\cite{clem08b}) are nothing but binary classification problems with the peculiarity that one has to predict exactly $K$ class 1 instances. Since a binary classification loss function is used for weak ranking problems, the notion of the angular breakdown point of \cite{zhao18} (resp. \cite[Def. 2+2']{zhao18} for kernel-based classification) is directly applicable, but the results of \cite{zhao18} are only valid if the loss function is a suitable surrogate of the $0/1-$loss function since continuity is assumed there. As for the outlier scheme, note that the number $K$ leads to an additional constraint in the proposed outlier set in the proof of \cite[Thm. 2]{zhao18}. Based on the mere classification loss, we can produce outliers that lead to a breakdown of the coefficient in terms of the OIBDP. The proofs are in the Appendix.

\begin{cor} \label{weakinflosscor} If the (classification) loss function satisfies Assumption \ref{infass}, the upper bound of the sample and population OIBDP for weak continuous ranking is $p/n$ for $p \le K-1$ and not existent otherwise. \end{cor}

\begin{thm} \label{weakrankthm} For the weak continuous ranking problem with the $0/1-$loss function, the OIBDP \\
\textbf{a)} is given by $m/n$ for $m=\lfloor K/2 \rfloor+1$ for $p=1$, \\
\textbf{b)} is bounded from above by $K/n$ for $1<p<K$, \\
\textbf{c)} does not exist for $p \ge K$.
\end{thm}

We abstain from detailing out possible results for localized binary and localized $d-$partite ranking problems as well as for weak binary ranking problems. The reason is that these problems are essentially ill-posed from the perspective of the OIBDP. The reason is that when localizing, the top $K$ instances may all have the same label which makes them indistinguishable and therefore not rankable in any sense. We suggest to focus only on the classification part, inevitably requiring to measure the robustness in terms of the angular BDP for binary (\cite{zhao18}) or of the angular BDP for $d-$partite localized ranking (\cite{qian}).

\section{Discussion} \label{discusssec}

\subsection{The non-existence issue}

\begin{ex} Let $p>1$ and assume that $X_{ij}=0 \ \forall j \ne j_0$ and $\beta_j=I(j=j_0)$ for some $1 \le j_0 \le p$. Then, it suffices to use the worst-case outlier configuration from Fig. \ref{outliers} only on the $j_0-$th axis. Note that we cannot guarantee that our estimated coefficient maintains the zero components but however, $\beta_j \hat \beta_j \le 0$ is clearly satisfied. \end{ex}

The computed breakdown points depend on the dimension $p$ and generally grow with $p$. Even worse, if $p$ is at least as large as $n$ resp. $K$, a breakdown can no longer be achieved. However, the tides turn once sparsity of the true underlying model is assumed as the example above showed. 

\begin{Def} The linear model $Y=X\beta$ is called sparse with true dimension $q$ if $||\beta||_0=q$. In this setting, denote the set of relevant variables by $S^0$. \end{Def}

If the outlier scheme exactly knows which $q$ predictors are relevant (we may call the outliers \textbf{''oracle outliers''} here), the outlier scheme is only applied to the corresponding $q$ axes. This is true since for $\beta_j=0$ for $(p-q)$ components, the sign of the perturbed coefficient $\hat \beta_j$ obviously satisfies $\beta_j \hat \beta_j \le 0$, so enforcing a sign-reversal of the $q$ true non-zero coefficients only indeed is a breakdown. 

\begin{cor} Let $n$ be fixed and let $q \le n-1$ resp. $q \le K-1$ be the true dimension of the linear model, i.e., $||\beta||_0=q$. Then the order-inversal breakdown point for every ranking problem that we considered in this work exists. \end{cor}

\begin{cor} Let $q=q(n)=b_nn$ such that $b_n \in ]0,1[ \ \forall n$. Then the asymptotic order-inversal breakdown point for all non-localized ranking problems considered in this work exist. For localized ranking problems with $K=K(n)$ with $K(n)/n \rightarrow d \in ]0,1]$, we have to assume that $b_n \rightarrow b<d$. \end{cor} 

We are aware of the fact that very high-dimensional true models for which $q \ge n$ holds cannot break down in the sense of the OIBDP. In many situations, one can reduce this dimension to $q'<n$ by only considering the most relevant predictors (e.g. \cite{bu10}), although there are situations in which more than $n$ selected predictors are desired (e.g. \cite{wang11}). We do not think of this issue as being a weakness of our OIBDP notion since the OIBDP is quite intuitive and since the global nature of ranking problems that take at least pairs of instances into account and no single instances defines a significantly different setting than for example regression for which higher dimensions generally reduce the BDP. The OIBDP can still be used to compare the robustness of competing algorithms by considering the $q<n$ case which identifies which algorithm is more robust. 

\begin{rem} \label{negrem} One could ask why one cannot just multiply the responses with $(-1)$ in order to achieve a breakdown which also holds for SVR-type ranking estimators below in Lemma \ref{svrlem}. Honestly speaking, from an algorithmic perspective, we believe that one can indeed let the ranking estimator break down for any reasonable algorithm using this outlier scheme, making the OIBDP indeed existent for any true dimension $q$ (and therefore, letting it also exist for non-sparse true models). However, from a theoretical perspective, there is no evidence that one cannot result in a non-broken coefficient since the solution set, i.e., the set of all coefficients that optimize a ranking loss for the data set with the negated responses, does not only consist of broken coefficients but also of non-broken ones. The argument is the same as in Thm. \ref{supbdpinf} that sign-inverting some but not all coefficient components may already lead to a perfect ranking prediction on the contaminated sample, so there is no guarantee that all components would be enforced to be sign-inverted. \end{rem}

\subsection{The sample OIBDP}

In principle, we could directly transfer the results from this paper that restricted themselves to the population OIBDPs to the sample OIBDP setting. However, the sample OIBDP heavily relies on the quality of the estimator $\hat \beta$ on the clean data set (note that the sample angular BDP from \cite{zhao18} faces the same issue). Let us first recapitulate an important definition that can be found for example in \cite{bu}. 

\begin{Def} \label{varselcons} Assume that $S^0$ is the true set of variables and $\hat S$ is the set of parameters selected by the model selection procedure. The model selection procedure is \textbf{variable selection consistent} if $ P(\hat S=S^0) \longrightarrow 1$  for $n \rightarrow \infty$. \end{Def} 

Variable selection consistency is a strict assumption. Theoretical results often cover a relaxed property, the so-called screening property (\cite{bu}) which indicates that the set of selected variables contains the set of true variables asymptotically, so variable selection consistency is the special case of equality of the two sets. Note that even sophisticated algorithms like $L_2-$Boosting (\cite{bu03}, \cite{bu06}) fail to be variable selection consistent (\cite{vogt}). %As for the Lasso, see \cite{bu} for the requirements for variable selection consistency.  \ \\

If the model selection procedure therefore has only the screening property, it essentially selects too many variables, making the ranking problem artificially more robust. We do not think that this is reasonable since this seeming robustness would result from the deficiencies of the applied model selection algorithm. Therefore, we only consider variable selection consistent model selection which enables the following asymptotic result.

\begin{cor}  Let $q=q(n)$ such that $q(n)$ with $q(n)/n \rightarrow b \in ]0,1[$ resp $q(n)/n \rightarrow b<d$ for $K(n)/n \rightarrow d \in ]0,1]$. Then the asymptotic sample order-inversal breakdown point for all non-localized ranking problems resp. localized ranking problems considered in this work exist provided that the estimated coeffient is computed using a variable selection consistent procedure. \end{cor}

Note that variable selection consistency is an asymptotic property, so for fixed $n$, even a procedure that satisfies this property can produce an estimated coefficient which selects irrelevant variables or which misses relevant variables. %Theoretically, one also could use a similar outlier scheme that is located at the axes which correspond to the non-zero entries of the estimated coefficient which however would require to know the estimated coefficient on the clean data which we do not think to be a reasonable assumption in contrast to an oracle property which has information about the true non-zero coefficients. On top of that, outliers can clearly affect the selected model. 
Summarizing, we believe that sample versions of BDPs are not very informative and should always be considered with caution.

\subsection{Lower bounds for the OIBDP}

Lower bounds for the OIBDP in the sense that one asks for example in the situation of Lem. \ref{infbdplemma} where we assumed that the original data points supply non-broken coefficients most (i.e., that the data are linearly rankable) if there is any lower OIBDP value that holds with high probability on real data (where linear rankability may not hold) cannot be computed universally due to numerous reasons. 

First, the original data contain some noise so even if they would follow some linear model with some true $\beta$, the observed response values would differ from the ideal response values so that linear inrankability can occur by chance. However, if one had a model $Y_i=s_{b,\beta}(X_i)+\epsilon_i$ for some stochastic error term $\epsilon_i$, the probability that linear inrankability occurs does not only depend on the error distribution but also on the $X_i$ and on $\beta$. For example for $p=1$, if $\beta>0$ is very large, the probability that the errors make the data points linearly inrankable would be smaller than for some smaller $\beta>0$ if the predictors are kept fixed. For a fixed $\beta$, the probability that linear inrankability occurs would also be lower if there are large distances between the predictors since the expected responses then would be better separated. 

We also already mentioned in Rem. \ref{gradrem} that the underlying numerical algorithm itself may affect the OIBDP. Due to these reasons, we think that if one had a concrete algorithm, a given data set and a good intuition of the error distribution and the true coefficient, one may would be able to compute lower bounds for the OIBDP, but evidently, there is no chance to provide universal results.

\subsection{Practical implications}

The results from this work indicate that bounded loss functions lead to more robust ranking problems than unbounded loss functions. This is not surprising and coincides with the well-known results from robust regression and robust classification where redescenders, i.e., loss functions whose gradient in absolute value redescends to zero so that the loss functions asymptotically grow until reaching a constant, are proposed. 

As for ranking losses, we always assumed that $\lim_{u \rightarrow \infty}(L(u))=0$. Therefore, bounded loss functions with $\lim_{u \rightarrow -\infty}(L(u))=C_l<\infty$ as we assumed in several theorems are in fact redescenders. The problem is that redescenders are non-convex which makes numerical optimization difficult. In fact, almost every existing ranking algorithm works with convex and therefore unbounded surrogate loss functions. Nevertheless, non-convex optimization has already been addressed in for example robust regression, so developing a robust ranking algorithm is definitely possible, although, due to the global nature of ranking loss functions, the computational complexity can be assumed to be very high. Therefore, providing a robust ranking algorithm is beyond the scope of this work. 

On the other hand, a standard robustification technique is trimming on which many successful machine learning algorithms are based, most prominently the LTS (\cite{rous84}) or the SLTS (\cite{alfons13}). These trimming techniques are based on the in-sample losses, i.e., one iteratively identifies the relative $(1-\alpha)-$fraction of instances with the lowest in-sample loss, updates the model by fitting it on these instances, checks again which instances provide the lowest loss and so forth. We want to point out why this trimming technique is not trivially applicable to ranking. 
\begin{figure}[H] \begin{center}  \includegraphics[width=4cm]{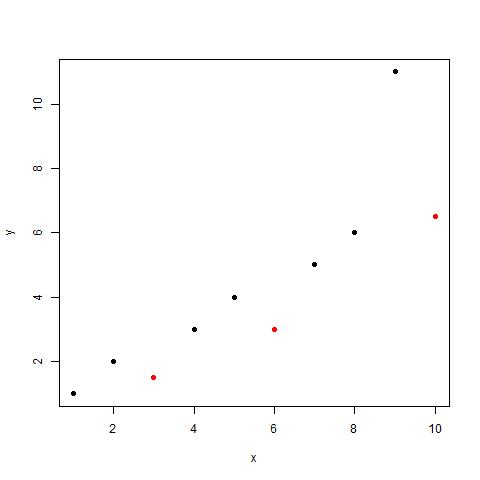}  \end{center}  \label{trim} \end{figure} 
Looking at Fig. \ref{trim}, there are three instances that contradict a positive ranking coefficient, colored in red. If one would apply trimming with a trimming rate of $\alpha=0.2$, the first question that arises is which of the three red points should be discarded since the indicator loss would make them indistinguishable in terms of the loss, so one would have to pick two of them randomly. Usually, one considers a surrogate of the indicator loss which would clearly discard the points (6,3) and (10,6) because the loss when comparing these instances with their left neighbors leads to the values \begin{center} $ \displaystyle (3-4)(6\beta-5\beta)=-\beta \ \ \ \text{resp.} \ \ \ (6-11)(10\beta-9\beta)=-5\beta$ \end{center} as input for the surrogate loss which is negative for all $\beta>0$, so due to the monotonicity assumption, these pairs lead to the highest losses. Since the comparison of the points (5,4) and (9,11) with their left neighbors does not produce a loss, one would learn that indeed the points (6,3) and (10,6) are problematic from the perspective of ranking. If the trimming rate would be $\alpha=0.3$, one would discard all three red points. 

Although this argumentation seems to be logical, it is in fact highly misleading. Regarding the points themselves, all of them except for (9,11) are likely to have been created by a linear model, so (9,11) would appear as a regression outlier. From the perspective of ranking, this outlier lets however the point (10,6) \textbf{appear as an outlier which can be interpreted as a swamping effect (see e.g. \cite{rous11})}. Therefore, in contrast to regression or classification problems where each instance can be treated individually and where the in-sample losses are instance-specific, the \textbf{globality of ranking prevents from applying trimming techniques in the usual way}. 

One possible remedy, although not very popular in the ranking community, is to use a plug-in approach, i.e., one applies a regression algorithm and uses the regression predictions for the ranking prediction, in other words, the regression function serves as scoring function. Approaches in this direction have been proposed by \cite{sculley} who however combines a regression and a ranking loss while \cite{mohan11} solely consider the squared loss. In this spirit, robust ranking may be achievable by robust regression, i.e., one could perform algorithms like SLTS and use its predictions for the ranking prediction. This will be an interesting topic for future work.

\section{Outlook: SVM-type approaches} \label{furthersec} 

A large class of ranking algorithms are from the SVM-type which potentially operate in infinite-dimensional reproducing kernel Hilbert spaces (RKHS). At the first glance, such methods would be problematic for a ranking BDP since even finite-dimensional RKHSs like the ones induced by polynomial kernels would seemingly be prone to hurt the condition $p<n-1$. The angular BDP from \cite{zhao18} has already been extended to kernel-based classification methods where they require the angle between the linear expansion (due to the representer theorem, e.g., \cite{schol01}) of the true function resp. the solution computed on the contaminated data set, measured by the norm in the corresponding RKHS, to be non-positive. Similarly, due to the component-wise nature of our OIBDP and the representer theorem, we can propose a reasonable definition of a BDP for kernel-based ranking estimators which is similar as the angular BDP from \cite[Def. 2+2']{zhao18}.

\begin{Def}[Order-inversal breakdown point for kernel-based ranking] \label{bdprankdefkernel} Assume that the true model has the form \begin{equation} \label{reprthm}  f(x)=\sum_{i=1}^n \sum_{j=1}^n (\alpha_{i,j}-\alpha_{i,j}^*)(K(X_i,x)-K(X_j,x))+b  \end{equation} for some kernel $K$, an intercept term $b$ with $|b|<\infty$ and coefficients $\alpha_{i,j}, \alpha_{i,j}^* \ge 0$.  \\
\textbf{a)} The \textbf{population order-inversal breakdown point for kernel-based ranking} is defined by \begin{center} $ \displaystyle \check \epsilon(f,Z_n):=\min\left\{\frac{m}{n} \ \bigg| \ \hat f(Z_n^m) \in S_{\cap}^-\right\}, \ \ \ S_{\cap}^-:=\bigcap_{k: f_k \ne 0}\{\tilde  f_k \ | \ \langle \tilde f_k, f_k \rangle_{\mathcal{H}} <0\} $ \end{center} where $\mathcal{H}$ is the RKHS corresponding to $K$ and where $f_k$ is the $k-$th component of $f$. \\
\textbf{b)} The \textbf{sample order-inversal breakdown point for kernel-based ranking} is defined by \begin{center} $ \displaystyle \check \epsilon(\hat f,Z_n):=\min\left\{\frac{m}{n} \ \bigg| \ \hat f(Z_n^m) \in \hat S_{\cap}^- \right\}, \ \ \ \hat S_{\cap}^-:=\bigcap_{k: \hat f_k \ne 0}\{\tilde  f_k \ | \ \langle \tilde f_k, \hat f_k \rangle_{\mathcal{H}}<0\} $ \end{center} \end{Def} 

Standard SVM classification solutions do not invoke the $\alpha^*-$coefficients. For ranking algorithms that solely invoke $\alpha-$coefficients, w.l.o.g. set $\alpha^*_{i,j}=0$ for all $i$ to consistently cover both cases with the definition of the OIBDP for kernel-based ranking. This general assumption covers the SVM-type ranking approaches like \cite{joachims02}, \cite{cao}, \cite{bref}, \cite{pahi} and \cite{tian} where for example \cite{herb} let the class label enter as factor and \cite{rak04} let the indices $i$ and $j$ run through all positive resp. negative instances. Additional constraints for the coefficients like upper bounds as considered in \cite{rak04} are not relevant in our BDP setting while particular index sets are covered by setting the coefficients of the remaining summands to zero. As for the $\alpha^*-$coefficients, many of the existing ranking algorithms are tailored to bipartite ranking and essentially approximate the conditional probability $\eta(x):=P(Y=1|X=x)$ which relates ranking problems and regression algorithm like support vector regression (SVR). 

There already exist sparse SVMs for ranking, see \cite{tian}, \cite{pahi10}, \cite{lai13} and \cite{laporte}, but there is no guarantee that the selected number of features would be smaller than $n$. However, considering SVR techniques, due to the requirement that the coefficients $\alpha_i$ and $\alpha_i^*$ have to be non-negative, we can conclude that for SVR-type algorithms, we need to enforce that $\sign(\alpha_i-\alpha_i^*)$ switches for every $i=1,...,n$ while preserving the sign of the differences of the features or kernelized features or vice versa. Studying the quantitative robustness of SVMs and SVRs in terms of the OIBDP which requires a thorough investigation of the corresponding dual problems for the $\alpha-$ (and $\alpha^*-$)coefficients would exceed the scope of this work. However, we can state an enlightening result regarding standard SVR. Note that the proposed outlier scheme is the same as in \cite{zhao18}. 

\begin{lem}\label{svrlem} If $(\hat \alpha, \hat \alpha^*)$ is the solution to the standard SVR problem (see. e.g. \cite{elstat}) \begin{equation*} \begin{gathered} \min_{\alpha, \alpha^*}\left(\epsilon \sum_i (\alpha_i^*+\alpha_i)-\sum_i y_i(\alpha_i^*-\alpha_i)+\frac{1}{2}\sum_i \sum_j (\alpha_i^*-\alpha_i)(\alpha_j^*-\alpha_j)\langle x_i, x_j \rangle \right) \\ 0 \le \alpha_i, \alpha_i^* \le C, \ \ \ \sum_i (\alpha_i^*-\alpha_i)=0, \ \ \ \alpha_i \alpha_i^*=0 \ \forall i \end{gathered} \end{equation*} for some cost parameter $C$ and the cutoff $\epsilon$ from the $\epsilon-$insensitive loss function, $(\hat \alpha^*, \hat \alpha)$ is the solution of the SVR problem on the data where the signs of all responses were switched. \end{lem}

\begin{bew} The statement is easily seen since $(\hat \alpha^*, \hat \alpha)$ obviously satisfies the constraints. The first sum of the objective does not change, also the third sum does not change by switching the sign of the two factors. The negation of the second sum due to the sign switch of the responses is compensated by the sign switch of the coefficient differences, so the value of the objective for the solution $(\hat \alpha^*, \hat \alpha)$ on the manipulated data is identical to the value of the objective of the solution $(\hat \alpha, \hat \alpha^*)$ on the clean data and since a minimum is attained, $(\hat \alpha^*, \hat \alpha)$ is optimal.\begin{flushright} $_\Box$ \end{flushright} \end{bew}

The proof is simple but the statement is of major importance since it already proves the astounding fact that \textbf{there is no ''blessing of dimensionality'' for support vector regression regarding our OIBDP} since the same statement is true when $\langle X_i, X_j \rangle$ is replaced by $K(X_i,X_j)$, so even infinite-dimensional feature spaces do not prevent the OIBDP from existing. This is no contradiction to Rem. \ref{negrem} since Lem. \ref{svrlem} is tailored to the special case of SVR, so the statement does not transfer to other machine learning algorithms. We will not extensively study all existing SVM-type ranking algorithms but we state the following for one of the most important and pioneering ranking algorithms.

\begin{cor} The OIBDP of the ranking SVM algorithm from \cite{herb}, \cite{herb99} always exists. \end{cor}

\begin{bew} Similarly as in the proof of Lemma \ref{svrlem}, let $\hat \alpha$ be a solution of the corresponding dual optimization problem which is given in \cite[Eq. (68)]{herb99}. The objective function invokes a double sum where factors $Y_iY_j$ appear (note that $Y_i \in \{\pm 1\}$ for all $i$ in their work). However, switching the sign of all responses will not affect the objective function and therefore keep the solution. Due to the linear expansion of the weights given in \cite[Eq. (69)]{herb99}, all summands that form the weights are sign-switched, so the whole weight coefficient is sign-reverted since the features stay untouched. The same is true if the kernelized SVMs which are computed by maximizing the objective function given in \cite[Eq. (75)]{herb99}, so by the analog expansion of the weights, implicitly given in \cite[Eq. (79)]{herb99}, a breakdown is achieved. \begin{flushright} $_\Box$ \end{flushright} \end{bew}

As for a general statement of the OIBDP for kernel-based ranking, we refer to the results from \cite[Thm. 3+Prop. 3+Prop. 4]{zhao18} who proved upper bounds for  their angular BDP for kernel-based classification if unbounded loss functions and unbounded kernels are considered. Although their angular BDP is not identical to our OIBDP, they essentially sign-revert every summand in the corresponding representer theorem expansion by keeping the predictor values but by switching the sign of the respective responses. Interestingly, they derive very similar results as we did for the linear ranking setting for unbounded loss functions, i.e., the upper bound for the BDP is given by $\tilde p/n$ if $\tilde p$ is the dimension of the RKHS induced by the kernel, so the same problems concerning BDPs greater than 0.5 or even non-existent BDPs occur here. As for unbounded RKHS's, the idea of \cite{zhao18} is to consider the effective dimension, i.e., the dimension of the finite-dimensional subspace of the RKHS in which the true scoring function $f$ can be represented. The resulting upper bound is then again given by this number divided by $n$. 

We postulate that the OIBDP for the SVM-type ranking algorithms always exists and takes a value lower than 1. This assumption is motivated by the results in \cite[Ch. 4]{zhao18} and by Lemma \ref{svrlem}. However, due to the huge variety of SVM-type ranking algorithms, we leave rigorous results about their OIBDPs, both regarding upper and possible lower bounds and for bounded resp. unbounded kernels, open for future research.

\section{Conclusion} 

We introduced the order-inversal breakdown point for ranking and argued why neither the classical regression breakdown point nor the angular breakdown point for classification are appropriate for this setting. We then systematically studied the breakdown points for different types of ranking problems that we carefully distinguished. Our contribution includes least favorable outlier configurations and corresponding characterizations of the OIBDP as well as sharp asymptotic upper bounds, respecting all types of ranking problems that are appropriate for this setting combined with the extreme cases of unbounded loss functions and non-continuous indicator loss functions. 

One could argue that our BDPs may not be reasonable since cases with asymptotic BDPs of 1 or even cases where the BDP does not even exist arise. However, these problems are directly related to the sparsity of the underlying true model. Since a sparsity assumption is always encouraged in high-dimensional settings, relatively mild conditions on the growing behaviour of the predictor dimension allow for an OIBDP smaller than 1. 

Our results imply that robust ranking can be achieved by optimizing (non-convex) redescending surrogate losses, but we leave the derivation of a concrete algorithm of this type as well as studying the plug-in approach based on robust regression open for future research. We also shortly discussed an extension of our OIBDP for ranking for the case of SVM-type scoring functions and proved the existence of this BDP, even for infinite-dimensional feature spaces.  

\renewcommand\refname{References}
\bibliography{Biblio}
\bibliographystyle{abbrvnat}

\appendix

\section{Characterization of the OIBDP}

We can prove an analog to \cite[Thm. 2]{zhao18}. In our work, the result is of lesser importance since we cannot conclude that the OIBDP always exists which \cite{zhao18} indeed can for their angular BDP. Assume that we have a bounded loss function, i.e., $L(u) \le C_l<\infty$. Define \begin{center} $ \displaystyle G_{\lambda,m}^u(\tilde \beta,Z_{n-m})=G_{\lambda,n}(\tilde \beta,Z_{n-m})+\frac{n(n-1)-(n-m)(n-m-1)}{n(n-1)}C_l, $ \end{center}  as an analog to \cite{zhao18} with $G_{\lambda_n}$ from Eq. \ref{objregrank}. The second summand indicates the upper loss bound achieved due to the outliers, i.e., both due to the pairwise comparisons of outliers as well as due to the pairwise comparisons of outliers and non-outliers.  

\begin{thm} \label{charbdp} Let $\beta \ne 0_p$ and let the loss function be decreasing with $\lim_{u \rightarrow \infty}(L(u))=0$ and $\lim_{u \rightarrow -\infty}(L(u))=C_l<\infty$. Then the estimator does not break down in the sense of the population resp. sample OIBDP for ranking if and only if \begin{equation} \label{glguineq} \min_{\tilde \beta_1 \in \Delta_{BL}^+}(G_{\lambda,n}^u(\tilde \beta_1,Z_{n-m}))<\min_{\tilde \beta_2 \in \Delta_{BL}^-}(G_{\lambda,n}(\tilde \beta_2,Z_{n-m})) \end{equation} for $S_{\cap}^-$ as in Def. \ref{bdprankdef}, $S_{\cap}^+=\R^p \setminus S_{\cap}^-$ and $\Delta_{BL}^+=\{(b,\beta) \ | \ \beta \in S_{\cap}^+, |b|<\infty\}$ and $\Delta_{BL}^-$ analogously.\end{thm}

\begin{bew} We argue along the same lines \cite{zhao18} with some modifications but for making the proof self-contained, we detail out the steps. We restrict ourselves to the population setting since in the sample setting, one just has to replace the original coefficient by the coefficient estimated on the original data. We define the following set of outliers: \begin{equation*} \begin{gathered} \check Z_0^m(X,Y,\beta):=\{(X_i^0,Y_i^0) \ | \ X_{ij}^0=X_{ij}+c_{ij} \ \forall j: \beta_j \ge 0, X_{ij}^0=X_{ij}-c_{ij} \ \forall j:\beta_j<0, \\  0<c_{ij} <c_{kj} \ \forall i<k \ \forall j, |X_{ij}^0|>\max_i(|X_{ij}|) \ \forall i, Y_i^0>Y_k^0 \ \forall i<k, \max(Y_i^0)<\min(Y_i)\} .\end{gathered} \end{equation*} Graphically, this set is easily understood and is depicted for an example in Fig. \ref{outliers} in the proof of Lemma \ref{infbdplemma} for $p=1$. By construction of the $X_i^0$, we proceed along the cone where $x \beta$ is increasing ensuring that the magnitude of each component $X_{ij}^0$ exceeds the magnitude of each $X_{ij}$ and that the $X_i^0$ are different. However, the responses are defined such that they are descending with $i$ while the original coefficient $\beta$ would lead to an ascending ordering, i.e., the ordering is reverted and on top of that, the ordering of each outlier compared with each original observation is reverted. 

Consequently, for any $\tilde \beta_2 \in \Delta_{BL}^-$, we have \begin{center} $ \displaystyle L_{\lambda,n}(\tilde \beta_2,\check Z_n)=\min_{Z_m^0}(L_n(\tilde \beta_2,\tilde Z_n))=G_{\lambda,n}(\tilde \beta_2,Z_{n-m}) $ \end{center} where $\check Z_n=Z_{n-m} \cup \check Z_m^0$ for $\check Z_m^0 \in \check Z_m^0(X,Y,\beta)$ and $\tilde Z_n=Z_{n-m} \cup Z_m^0$ for any outlier set $Z_m^0$ since any such $\tilde \beta_2$ reverts the ordering on the original data but makes perfect predictions for all pairs of outliers and all pairs with one outlier and one original response. We cannot guarantee that any $\tilde \beta_1 \in \Delta_{BL}^+$ achieves the worst-case loss for the components $H_n$ and $F_n$, but by construction, for such a given $\tilde \beta_1$, there definitely exists an outlier set $\check Z_m^0$ such that \begin{center} $ \displaystyle L_{\lambda,n}(\tilde \beta_1,\check Z_n)=\sup_{Z_m^0}(L_n(\tilde \beta_1,\tilde Z_n))=G_{\lambda,n}^u(\tilde \beta_1,Z_{n-m}). $ \end{center} Although we cannot guarantee that a worst-case outlier set exists such that every coefficient that does not satisfy the breakdown criterion suffers supremal loss, for now we can only conclude that the estimator does not breakdown if \begin{center} $ \displaystyle \min_{\tilde \beta_1 \in \Delta_{BL}^+}(G_{\lambda,n}(\tilde \beta_1,Z_{n-m}))<\min_{\tilde \beta_2 \in \Delta_{BL}^-}(G_{\lambda,n}(\tilde \beta_2,Z_{n-m})). $ \end{center} Now, we are ready to prove the stated equivalence.  \\
\textbf{i)} Assume that the estimator does not break down, i.e., the computed estimator $\hat \beta_{\lambda}(\tilde Z_n)$ is contained in $\Delta_{BL}^+$ for any outlier set $Z_m^0$. Then, due to the fact that for any $\tilde \beta_1 \in \Delta_{BL}^+$, there exist an outlier set such that $\tilde \beta_1$ suffers the maximal loss $G_{\lambda,n}^u(\tilde \beta_1,Z_{m-n})$, it follows that the inequality in Eq. \ref{glguineq} indeed holds.  \\
\textbf{ii)} Assume that inequality \ref{glguineq} holds. The property $L_{\lambda,n}(\tilde \beta_1,\check Z_n) \le G_{\lambda,n}^u(\tilde \beta_1,Z_{n-m})$ obviously holds. The statement $L_{\lambda,n}(\tilde \beta_2,Z_n) \ge L_{\lambda,n}(\tilde \beta_2,\check Z_n)$ for $\tilde Z_n=Z_{n-m} \cup Z_m^0$ holds by construction for any outlier set $Z_m^0$ still holds so that we conclude \begin{center} $ \displaystyle \min_{\tilde \beta_1 \in \Delta_{BL}^+}(L_{\lambda_n}(\tilde \beta,\check Z_n)) \le \min_{\tilde \beta_1 \in \Delta_{BL}^+}(G_{\lambda_n}(\tilde \beta_1,Z_{n-m}))<\min_{\tilde \beta_2 \in \Delta_{BL}^-}(G_{\lambda,n}(\check \beta_2,Z_{n-m})) \le \min_{\tilde \beta_2 \in \Delta_{BL}^-}(L_{\lambda,n}(\tilde \beta_2,\tilde Z_n)) $ \end{center} where the strict inequality holds by assumption and where the last inequality holds since any outlier set from the worst case outlier set $\check Z_0^m(X,Y,\beta)$ lets no pair of an outlier and an original response suffer any loss for any coefficient from $\Delta_{BL}^-$ which is not guaranteed by general outlier sets. Therefore, the estimator does not break down since a coefficient from $\Delta_{BL}^+$ will be optimal, i.e., achieve the minimum loss. \begin{flushright} $_\Box$ \end{flushright} \end{bew} 

%\begin{rem} The main difference between the setting of \cite{zhao18} is that no worst-case outlier set exists such that any $\tilde \beta_1 \in \Delta_{BL}^+$ suffers the worst-case loss. However, an estimator does not break down if and only if for \textbf{any} outlier configuration of the respective cardinality it can be guaranteed that there exists such a $\tilde \beta_1$ which achieves a lower loss than any $\tilde \beta_2 \in \Delta_{BL}^-$ which is the reason why our statement follows. \end{rem} 

\section{Additional proofs and examples}

\subsection{Hard ranking}

The following example has been announced in Rem. \ref{outlierrem}.

\begin{ex} \label{outlierdice} We already thought of an outlier scheme with more dependencies where, here wlog. $p=2$, we do not use $X^{(3)}=(X_1'+2,X_2')$ but $X^{(3)}=(X_1'+1,X_2'+1)$, i.e., the outliers form the edges of a square. Then these four outliers already enable two comparisons along both axes which required five outliers using the proposed outlier scheme in the proof of Thm. \ref{supbdpind}. The next iteration would be to build a square with vertex length 2 where the outliers define the edges, the mid-points of the vertices and the middle point of the square which leads to nine comparisons along each axis. In general, we would have $k^2$ outliers and $k^2(k-1)/2$ axis-wise comparisons. The general strategy would construct a $p-$dimensional hypercube grid with $k^p$ outliers, leading to $k^p(k-1)/2$ axis-wise comparisons. This strategy can be beneficial for small $n$, for example in the case $p=3$ and $k=2$, we have 4 comparisons along each axis using only 8 outliers while our strategy before would require 12 outliers to beat this (9 outliers would only lead to 3 comparisons along each axis). However, it is easily revealed that this strategy does not work for larger $n$, for example in the case $p=2$ and $k=5$, we had 25 outliers and 50 comparisons, but with the strategy before, 24 outliers would already enable 66 comparisons. Maybe future research is able to nevertheless reveal a better strategy than ours proposed in the proof of Thm. \ref{supbdpind}.   \end{ex}

\begin{bew}[Proof of Cor. \ref{supbdpindcor}] \textbf{i)} We prove the statement by setting $k:=cn$ for some $c \in ]0,1]$. The condition for a breakdown is then \begin{equation} \label{compeq} \begin{gathered} \frac{cn(cn+1)}{2} \overset{!}{>} \frac{(n-pcn-1)(n-pcn-2)}{2} \\ \Longleftrightarrow c^2n^2+cn \overset{!}{>} n^2-2pcn^2-3n+p^2c^2n^2+3pcn+2 \\  \Longleftrightarrow n^2(c^2+2pc-p^2c^2-1) \overset{!}{>} n(3pc-c-3)+2 \overset{n>0}{\Longleftrightarrow} n(c^2+2pc-p^2c^2) \overset{!}{>} 3pc-c-3+\frac{2}{n} \end{gathered}  \end{equation} and since $3pc-c-3$ is fixed and $2/n \rightarrow 0$ asymptotically, we only have to guarantee that the bracket is positive, i.e., \begin{center} $ \displaystyle c^2+2pc-p^2c^2>0 \Longrightarrow c=\begin{cases} \frac{-1-p}{1-p^2} \\ \frac{1-p}{1-p^2}=\frac{1}{p+1} \end{cases} $ \end{center} where the first case is not meaningful since it contradicts $c>0$. Therefore, we asymptotically set $c^*=1/(p+1)$ and get \begin{center} $ \displaystyle m^*=1+c^*np=1+\frac{p}{n(1+p)} \Longleftrightarrow \frac{m^*}{n}=\frac{1}{n}+\frac{p}{1+p} $ \end{center} which asymptotically equals $p/(p+1)$ as stated. 

\textbf{ii)-iii)} The second statement is obvious since for any such sequence $(b_n)_n$ eventually leads to a diverging number $p(n)$ of variables with $p(n)/(1+p(n)) \rightarrow 1$. Note that in the last statement in the last line in Eq. \ref{compeq}, the right hand side cannot grow indefinitely with $p$ since $pc<1$ (due to the requirement that $n-pcn-2 \ge 0$), so the same computations as for static $p$ hold. The third statement already has been discussed for the fixed $p \ge n$. \begin{flushright} $_\Box$ \end{flushright} \end{bew}

The following example provides finite-sample upper bounds for the OIBDP for hard ranking with the indicator loss in selected scenarios.

\begin{ex} We simulate the BDP for different $p$ and for a sequence of values for $n>p$: \begin{figure}[H] \begin{center} \includegraphics[width=4cm]{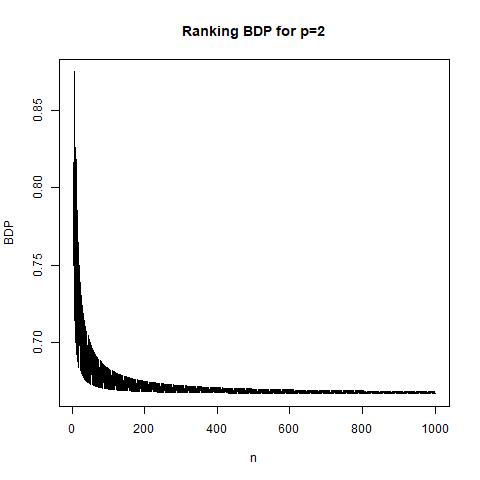} \includegraphics[width=4cm]{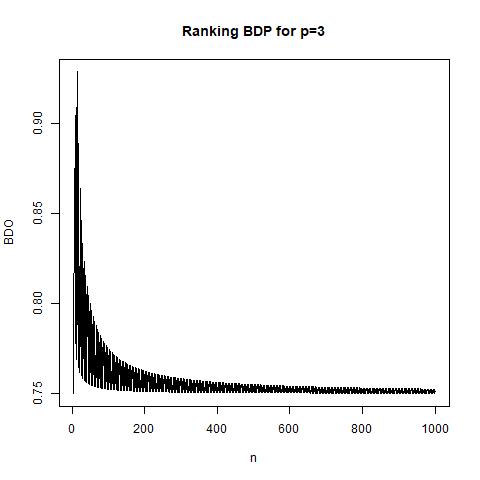} \includegraphics[width=4cm]{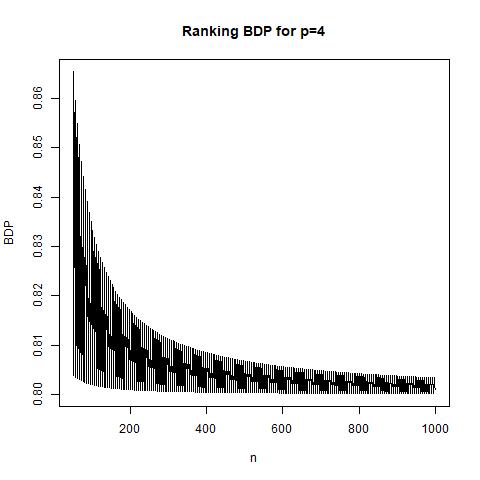} \includegraphics[width=4cm]{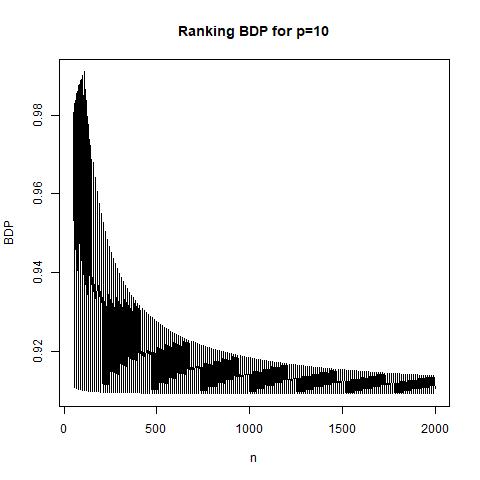} \caption[Upper bounds for the OIBDP]{Plots of upper bound for the BDPs for ranking for different $p$} \label{checkmplot} \end{center} \end{figure} \end{ex}

Let us outline an artificial case where the outlier flexibility is severely hindered so that the expected OIBDP is better suited.

\begin{ex} \label{compex} Note that we implicitly assumed \textbf{open} regressor resp. responses spaces when constructing the outlier set. If we have compact regressor or response sets, there is no guarantee that the outlier scheme is applicable. %since we require that (wlog. $\beta_j>0$ for all $j$) the outlier responses are lower than all original response which is impossible if there are too many original responses located at the minimal response value. In fact, if the number of such responses would be smaller than the respective $m^*$ in the theorems, one would just replace these points by the outliers. However, a general remedy would be more difficult. Since we forbid the case that all original responses are equal (since either $\beta<0$ or $\beta>0$ for the true coefficient), there are at most $(n-1)$ original responses located at the minimal value. This case could be easily handled since we can mirror our outlier scheme such that the outlier responses are higher than all original responses while the predictor values are lower. The most difficult case however arises if $\lfloor n/2 \rfloor$ of the original responses are located at the minimal response value and $\lfloor n/2 \rfloor$ at the maximal response value and if additionally the predictor values are located at or near the allowed maximal resp. minimal values. \ \\
To illustrate this setting for $p=1$, let again $\beta>0$ and let the original data be linearly rankable according to this coefficient. More precisely, assume the (very artificial) case that $n$ is even and $n/2$ points are given by $(\max(\mathcal{X}),\max(\mathcal{Y}))$ and the other half of the points at the respective minima. Assuming a bounded loss function, wlog. the indicator loss function, we have no choice but to replace one of these clusters completely by keeping the regressor value but by moving the response to the other extremum of the response space. Now, the usual outlier scheme that we already introduced is no longer applicable. The only chance we have is to start by replacing wlog. the whole upper cluster by $n/2$ outliers according to the scheme $(X_i^0,Y_i^0)=(\max(\mathcal{X})-\epsilon_i,\min(\mathcal{Y})+\epsilon_i)$ for $\epsilon_1>\epsilon_2>...>\epsilon_k>0$ with $\epsilon_1$ being small enough to let the first outlier be contained in the open interior of $\mathcal{X} \times \mathcal{Y}$. This strategy is depicted in Fig. \ref{outliercomp} where we jittered the points at the left corner only to make them visible. 
\begin{figure}[H] \begin{center}  \includegraphics[width=4cm]{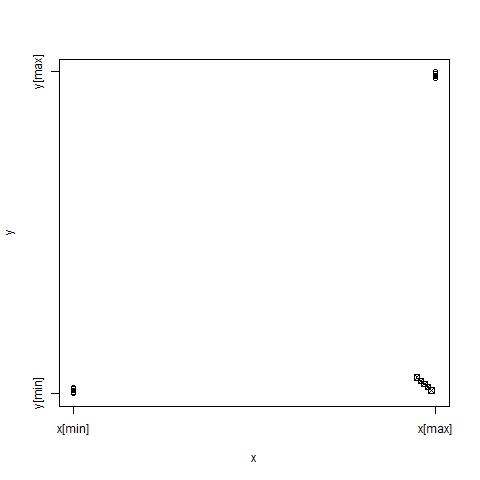}  \end{center} \caption[Worst-case outliers for $p=1$]{Worst-case outliers for $p=1$} \label{outliercomp} \end{figure} 
Tedious algebra reveals an asymptotic BDP of $p^2/(p^2+0.5)$. \end{ex}

\subsection{Hard binary and hard $d-$partite ranking problems}  

\begin{bew}[Proof of Thm. \ref{supbdpindbinary}] Let us illustrate the proof for $p=2$. Similarly as in the univariate case in Lem. \ref{infbdplemmabinary}, it does not suffice to generate axis-wise outliers along one direction (i.e., either for very large or very small $X_{\cdot,j}$ for axis $j$) but one has to generate outliers on both sides. More precisely, again assuming a starting point $X^*$, for each $k$ one has to produce one axis-wise outlier on axis $j$ where the $j-$th variable is greater than the $j-$th variable for all other data and where the response is -1 (wlog. let again $\beta_j>0$ for all $j$) and one outlier where the $j-$th variable is lower than the $j-$th variable for all other data with response 1. This leads to $m=1+2pk$ outliers per iteration. On each axis, there are $k$ outliers on each side, leading to $k(k+1)$ misrankings since the starting point $X^*$ either has response 1 or -1, leading to $kl$ additional misrankings. In contrast, the sign-reverted coefficient potentially causes misrankings between all remaining $(n-2pk-1)$ original data points, so the formula \ref{kstarhardbinary} is proven. 

Clearly, if no remaining data points would be available, an early stopping strategy is applicable, i.e., it would suffice to let the starting point have response 1 and to only consider one axis-wise outlier with larger regressor value and response -1, so for $p \le n-1$, the BDP always exists. 

Again, since there essentially is no difference in the robustness of bipartite and $d-$partite ranking problems as already discussed in Rem. 6.1 and the proof of Lemma \ref{bipkpartrem}, the results directly transfer to $d-$partite ranking problems. \begin{flushright} $_\Box$ \end{flushright} \end{bew}

\begin{bew}[Proof of Cor. \ref{supbdpindcorbin}] Statement i) follows using some algebra as in similar statements before, ii) is true since the coefficient in i) converges to 1 for growing $p$ and iii) has already been discussed. \begin{flushright} $_\Box$ \end{flushright} \end{bew}

\subsection{Localized ranking problems} 

\begin{bew}[Proof of Cor. \ref{locinflosscor}] The ranking part directly follows the proof of Thm. \ref{supbdpinf}. As for the case of an unbounded classification loss function, we propose a similar construction but only need one of the clean instances of class -1 as starting point for it. Since this instance is of class -1 and the outliers are constructed to be of class 1 but only differ from the starting point by the value of one component, the respective coefficient has to be sign-switched in order to let the respective outlier be classified as class-1-instance (otherwise, the classification loss could be arbitrarily high by letting the outlier response diverge). In contrast to the ranking part where only the top $K$ instances are compared so that the starting point already has to be a class-1-instance, we can generate $K$ axis-wise outliers in the classification part, so the BDP exists unless $p>K$. \begin{flushright} $_\Box$ \end{flushright} \end{bew}

\begin{bew}[Proof of Cor. \ref{supbdpindloccor}] Statement i) is trivial, statement iv) has already been discussed (since the BDP converges to 1 for $d \rightarrow 1$, there is no second regime-switching point as in the univariate case), the formulae in ii) can be easily computed and iii) follows directly. As for $d_0$, see Fig. \ref{asylocbdp2} for illustration where the black curve corresponds to the first formula in ii) and the red curve to the second formula. \begin{figure}[H] \begin{center}  \includegraphics[width=4cm]{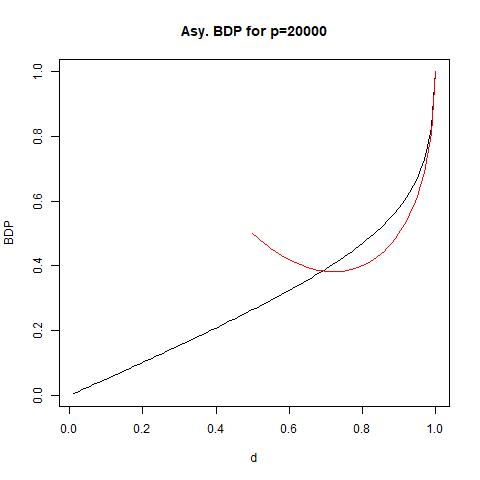}  \end{center} \caption[Asymptotic BDP for $K \le n/2$ and $K>n/2$]{Asymptotic BDP for $K \le n/2$ and $K>n/2$} \label{asylocbdp2} \end{figure} \begin{flushright} $_\Box$ \end{flushright} \end{bew}

\subsubsection{Localized continuous ranking problems on $\widehat{Best_K}$} 

Alternatively, we can localize the ranking loss on $\widehat{Best_K}$, i.e., the goal is to ensure that the instances that are predicted to be at the top of the list are ranked correctly, although these instances may not be the true top instances. We provide similar results. 

\begin{lem} \label{infbdplemmalochat}  Let $p=1$. For the localized continuous ranking problem with 0/1-loss for classification and the indicator loss function for ranking where the latter is based on the predicted best instances with indices in $\widehat{Best_K}$, the sample and population OIBDP for ranking  \\
\textbf{i)} is given by \begin{equation} \label{bdpcondlocbesthat} \frac{\check m}{n}, \ \ \ \check m=\min\left(K,\min\left\{k \ \bigg| \ \frac{n-K}{n} \cdot \frac{2(K-m)}{n}+\frac{(K-m)(K-m-1)}{2n(n-1)}<\frac{n-K}{n} \cdot \frac{2m}{n} \right\}\right)  \end{equation} for $K \le (n+m)/2$,  \\
\textbf{ii)}  is given by \begin{equation} \label{bdpcondlocbest2hat}  \begin{gathered} \frac{\check m}{n}, \ \ \ \check m=\min\{k \ \bigg| \ \frac{n-K}{n} \cdot \frac{2(n-K)}{n}+\frac{(K-m)(K-m-1)}{2n(n-1)} \\ <\frac{n-K}{n} \cdot \frac{2m}{n}+\frac{1}{n(n-1)}\left[\frac{m(m-1)}{2}+m(K-m)\right] \}   \end{gathered} \end{equation} for $K \in [(n+m)/2,n-m]$,  \\
\textbf{iii)} is given by Eq. \ref{bdpcond} in Lemma \ref{infbdplemma} where $n$ in the definition of $\check m$ is replaced by $K$ for $K \ge n-m$. \end{lem}

\begin{bew} The situation here is inherently different from the case that the ranking performance is computed on $Best_K$. We have to distinguish carefully between the two outlier schemes in Fig. \ref{outliers} and Fig. \ref{outliersloc}. 

Let us start with the case that $K<n/2$. The misclassification rate is obviously not affected by localizing the ranking performance on $\widehat{Best_K}$, so the formulae from Lemma \ref{infbdplemmaloc} remain valid. As for the misrankings, if we consider the outlier scheme as in Fig. \ref{outliersloc}, we will not produce any misranking for the original coefficient. This is true since for $\beta>0$, the rightmost instances are predicted to be the best ones and the ordering of their responses is correctly predicted as ascending. In contrast, any negative coefficient produces a complete inversion of the ordering of the remaining $(K-m)$ original instances, i.e., $(K-m)(K-m-1)/n$ misrankings. Let us now consider the outlier scheme as in Fig. \ref{outliers}. First note that the number of necessary outliers to produce a breakdown cannot exceed $K$ since we can achieve a zero loss for the broken coefficient using the outlier scheme in Fig. \ref{outliersloc} while the loss for the coefficient of the original sign is greater than zero due to the classification part. Now, the outliers according to Fig. \ref{outliers} cause $m(m-1)/2$ misrankings on the $m$ rightmost instances and additional $m(K-m)$ misrankings for any pairs of an outlier and one of the remaining rightmost $(K-m)$ instances. Any negative coefficient however again produces $(K-m)(K-m-1)/2$ misrankings on the intermediate $(K-m)$ instances. On the other hand, while the original coefficient only makes $m$ misclassifications, the broken coefficient misclassifies the maximum number of $K$ instances. 

Now, we have to argue which of the proposed outlier schemes applies. The answer is that it depends on $K$. Still assuming $K \le (m+n)/2$, we can observe that the classification loss is constant for the outlier scheme from Fig. \ref{outliers} for any $\beta<0$ for only a small additional ranking loss for $\beta>0$. We deduct that for $K \le (n+m)/2$, one should use the outlier scheme as in Fig. \ref{outliersloc} (asymptotically, it can be shown by numerical evaluation that the required number of outliers is always larger for the outlier scheme as in Fig. \ref{outliers} for $K=dn$ for $d \le d_0 \approx 0.692291$) and for $K>(n+m)/2$, we should use the outlier scheme as in Fig. \ref{outliers}. This proves formula \ref{bdpcondlocbesthat} and part i) as well as formula \ref{bdpcondlocbest2hat} and part ii) where the dependence on $K$ has already been discussed in the proof of Lemma \ref{infbdplemmaloc}.  

Finally, note that once $K \ge m-n$, the broken coefficient will only misclassify $(n-K)$ instead of $m$ instances, so the misclassification loss is equal for $\beta>0$ and $\beta<0$. Since only the ranking part remains which is the same as in the hard ranking setting with $n$ replaced by $K$, statement iii) is valid.  \begin{flushright} $_\Box$ \end{flushright} \end{bew}

\begin{cor} \textbf{i)} Asymptotically, a fixed $K$ would lead to a BDP of zero. For $d=1$, we get the asymptotic BDP of $1-\sqrt{0.5}$ as for hard ranking. \\
\textbf{ii)} For $K=K(n):=dn$ for $d \in ]0,d_0]$ for $d_0 \approx 0.5774659$, we can conclude that for $m=cn$, we have \begin{center} $ \displaystyle c^*=4-3d-\sqrt{16-28d+12d^2}$ \end{center} which takes values in $]0,0.30993]$ and is strictly monotonically increases with $d$. \\
\textbf{iii)} For $K=K(n):=dn$ for $d \in [d_0,d_1]$ with $d_1 \approx 0.773455$, we can conclude that for $m=cn$, we have \begin{center} $ \displaystyle c^*=1-\sqrt{-1+4d-5d^2/2}$ \end{center} which takes values in $[0.2265413,0.30993]$ and has its minimum at $d_1$. \\
\textbf{iv)} For $K=K(n):=dn$ for $d \in [d_1,1]$, we have the asymptotic BDP $d(1-\sqrt{0.5})$.  \end{cor} 

\begin{bew} Along the same lines as the proof of the Corollary \ref{infbdplemmaloccor}. The asymptotic BDP in dependence of $d$ is depicted in Fig. \ref{asylocbdphat}. \begin{figure}[H] \begin{center}  \includegraphics[width=4cm]{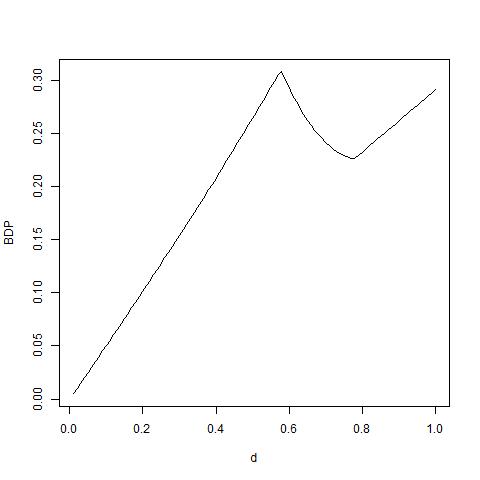}  \end{center} \caption[Asymptotic BDP for $K \le n/2$ and $K>n/2$]{Asymptotic BDP for $K \le n/2$ and $K>n/2$} \label{asylocbdphat} \end{figure} \begin{flushright} $_\Box$ \end{flushright} \end{bew}

\begin{thm} \label{supbdpindloc} Let wlog. the indicator loss functions for both the classification and the ranking part be used and let $p \ge 2$. Then, the upper bound for the OIBDP for localized ranking, localized on $\widehat{Best_K}$, is given by \begin{equation} \label{kstarlochat}   \begin{gathered} \frac{m^*}{n}, \ \ \ m^*=\min(1+pk^*,K), \\ k^*=\min\left\{k \ \bigg| \  \frac{n-K}{n}\frac{2(K-pk-1)}{n}+\frac{(K-1-pk)(K-2-pk)}{2n(n-1)}<\frac{n-K}{n}\frac{2(1+pk)}{n} \right\}  \end{gathered} \end{equation} for the case $K \le (n+m)/2$. For $K>(n+m)/2$, we have $m^*=1+pk^*$ where \begin{equation} \label{kstarlochat2}   \begin{gathered} k^*=\min\left\{k \ \bigg| \  \frac{n-K}{n}\frac{2(n-K)}{n}+\frac{(K-1-pk)(K-2-pk)}{2n(n-1)}<\frac{n-K}{n}\frac{2(1+pk)}{n}+\frac{k(k+1)}{2n(n-1)} \right\}  \end{gathered} \end{equation}  and for $K>n-m$, we get the same $k^*$ as in Eq. 5.3 in Thm. 5.2. This quantity always exists for $p \le K-1$. \end{thm}

\begin{bew} Along the same lines as the proofs of Thm. 5.1 and Lemma \ref{infbdplemmalochat} .  \begin{flushright} $_\Box$ \end{flushright} \end{bew}

\begin{cor} For the localized ranking problem where the localized ranking loss is computed on $\widehat{Best_K}$ and where both the classification and the ranking loss are indicator functions, we asymptotically conclude that \\
\textbf{i)} the BDP is zero for $K$ and $p$ being fixed, \\
\textbf{ii)} the BDP for $c \le 0.5$ tends to \begin{center} $ \displaystyle pc^*, \ \ \ c^*=\frac{4-3d}{p}-\sqrt{\frac{16-28d+12d^2}{p^2}} $ \end{center} and for $c>0.5$, it tends to the same quantity provided that $c \le 2(d-0.5)$, to \begin{center} $ \displaystyle pc^*, \ \ \ c^*=\frac{2p-pd}{p^2-1}-\sqrt{\frac{4p^2d-4p^2d^2-8d+5d^2+4}{(p^2-1)^2}} $ \end{center} for $c \in [2(d-0.5),1]$. For $p \rightarrow \infty$, the break-even point is given by $d_0 \approx 0.6923$. \\
\textbf{iii)} the BDP tends to the asymptotic BDP for hard ranking, i.e., $p/(p+1)$, for $d=1$. \\
\textbf{iv)} the BDP does not exist for $p=p(n)=b_nn$ with $b_n \rightarrow b \ge d$. \end{cor}

\begin{bew} Similar as in Cor. 7.3. Note that the first formula in ii) tends to the respective one in Cor. 7.3 for $p \rightarrow \infty$ while the second formulae are already equal. \begin{flushright} $_\Box$ \end{flushright} \end{bew}

\subsection{Other ranking problems} 

\begin{bew}[Proof of Cor. \ref{weakinflosscor}] Follows the same argumentation as the corresponding corollaries before. The starting point can be an arbitrary instance from the bottom of the list, so the axis-wise outliers (with a response tending to infinity) enforce the respective coefficients to switch their sign in order to keep the classification loss low. \begin{flushright} $_\Box$ \end{flushright} \end{bew}

\begin{bew}[Proof of Thm. \ref{weakrankthm}] Statement a) is trivial since the number of misclassifications is then lower for the coefficient with the opposite sign than for the original coefficient. Statement b) is only a coarse bound which cannot be tightened due to the missing ranking loss part in the weak ranking problems. The only opportunity that we have is to use an instance from the bottom of the list and to place axis-wise outliers with sufficiently large responses so that they are at the top of the list around this starting point. In the worst case, we have to start from an original instance from the bottom of the list and replace all $K$ true top instances with such axis-wise outliers. Statement c) is obvious. \begin{flushright} $_\Box$ \end{flushright} \end{bew}

\end{document}